\documentclass[11pt]{article}
\usepackage{amsmath,amsthm,amssymb,amsfonts,fullpage,verbatim,bm,graphicx,enumerate,epstopdf,lscape,subfigure,multirow}
\usepackage{algorithm,algorithmicx}
\usepackage{algpseudocode}
\usepackage{boxedminipage,longtable}
\usepackage[dvipsnames,usenames]{color}
\usepackage[pdftex,colorlinks=true,pdfpagemode=none,urlcolor=blue,linkcolor=blue,
   citecolor=blue,pdfstartview=FitH,plainpages=true]{hyperref}
\usepackage{tikz,verbatim,bibentry}
\usepackage{bm,mathrsfs,mathabx}
\usepackage[comma,authoryear]{natbib}
\usepackage{thm-restate}

\newtheorem{lemma}{Lemma}
\newtheorem{theorem}{Theorem}
\newtheorem{proposition}{Proposition}
\newtheorem{corollary}{Corollary}

\theoremstyle{definition}
\newtheorem{remark}{Remark}

\newtheorem{assumption}{Assumption}
\usepackage{etoolbox}
\AtEndEnvironment{remark}{\hfill$\blacksquare$}

\def\abs{\mathrm {abs}}

\def\F{\mathcal F}

\def\N{\mathcal N}

\def\I{\mathrm{I}}
\def\II{\mathrm{II}}

\def\L{\mathcal L}

\def\ZZ{\mathbb{Z}}
\def\EE{\mathbb E}
\def\PP{\mathbb P}
\def\log{\mathrm {log}}
\def\RR{\mathbb{R}}
\def\d{\mathrm{d}}
\def\NN{\mathbb{N}}
\def\1{\textbf{1}}

\def\exp{\mathrm{exp} }

\def\argmin{\mathrm{argmin} }

\begin{document}

\title{\textbf{Estimation of High-dimensional Nonlinear Vector Autoregressive Models}  }

\author{Yuefeng Han, Likai Chen and Wei Biao Wu\footnote{
Yuefeng Han is Assistant Professor, Department of Applied and Computational Mathematics and Statistics, University of Notre Dame, Notre Dame, IN 46556. Email: yuefeng.han@nd.edu. Likai Chen is Assistant Professor, Department of Statistics and Data Science, Washington University in St. Louis, St. Louis, MO 63130. E-mail: likai.chen@wustl.edu. Weibiao Wu is Professor, Department of Statistics, The University of Chicago, Chicago, IL 60637. E-mail: wbwu@uchicago.edu. Han was supported in part by National Science Foundation grant DMS-2412578. Chen was supported in part by National Science Foundation grant EF-2222403 and DMS-2311251. Wu was supported in part by National Science Foundation grants DMS-2311249 and DMS-2027723.
}}
\date{University of Notre Dame, Washington University in St. Louis and The University of Chicago}


\maketitle

\begin{abstract}
High-dimensional vector autoregressive (VAR) models have numerous applications in fields such as econometrics, biology, climatology, among others. While prior research has mainly focused on linear VAR models, these approaches can be restrictive in practice. To address this, we introduce a high-dimensional non-parametric sparse additive model, providing a more flexible framework. Our method employs basis expansions to construct high-dimensional nonlinear VAR models. We derive convergence rates and model selection consistency for least squared estimators, considering dependence measures of the processes, error moment conditions, sparsity, and basis expansions. Our theory significantly extends prior linear VAR models by incorporating both non-Gaussianity and non-linearity. As a key contribution, we derive sharp Bernstein-type inequalities for tail probabilities in both non-sub-Gaussian linear and nonlinear VAR processes, which match the classical Bernstein inequality for independent random variables. Additionally, we present numerical experiments that support our theoretical findings and demonstrate the advantages of the nonlinear VAR model for a gene expression time series dataset.

\ \\
\noindent%
\end{abstract}

\section{Introduction} \label{section:introduction}

The increasing variety of scientific applications has created a growing need for employing a large set of time series (variables) to model complex social and physical systems. This demand arises from various fields, including genomics \citep{sharon2013time}, neuroscience \citep{moller2001instantaneous, pereda2005nonlinear,kato2006statistical}, social networks \citep{ait2015modeling}, economics \citep{barigozzi2017network}, environmental studies \citep{lichstein2002spatial}, and communication engineering \citep{baddour2005autoregressive}. For example, economic policymakers rely on large-scale models of economic indicators \citep{sims1980, bernanke2005measuring, banbura2010large}, as empirical evidence has shown that they improve forecasts and provide better estimates of how current economic shocks will propagate, which guides policy actions more effectively. Similarly, in genomics and neuroscience, the advent of high-throughput technologies has enabled researchers to collect measurements on hundreds of genes or brain regions \citep{shojaie2010discovering, seth2015granger}, facilitating comprehensive modeling and deeper insights into biological mechanisms. In social sciences, many key variables are not directly observable but can be inferred through related time series variables, enabling a more nuanced understanding of policy decisions \citep{lin2020regularized}. Given the wide availability of high-dimensional time series data, understanding their underlying dynamic patterns is crucial for improving practical applications in these domains.


A widely used and informative model for capturing linear temporal dependencies between time series is the vector autoregression (VAR) model. Properties of VAR have been extensively studied in low-dimensional settings; see \cite{lutkepohl2005new} for a comprehensive overview. Over the past decade, a growing body of literature has leveraged structured sparsity and regularized estimation frameworks to achieve consistent estimation of VAR parameters in high-dimensional settings. \cite{basu2015} investigated the theoretical properties of Lasso-penalized high-dimensional VAR models for Gaussian processes. Their result was extended to multi-block VAR models by \cite{lin2017} and to factor-augmented VAR models by \cite{lin2020regularized}. \cite{guo2016} introduced a class of VAR models with banded coefficient matrices, which was further developed into spatio-temporal VAR models by \cite{gao2019}. \cite{basu2019low} explored high-dimensional VAR models involving low-rank and group-sparse components in network structures. \cite{hall2018} studied regularized high-dimensional autoregressive generalized linear models, focusing on Bernoulli and Poisson distributions. Additionally, \cite{ghosh2019, ghosh2021strong} developed Bayesian VAR models and analyzed their posterior and strong selection consistency. For further related work, see \cite{zheng2019testing, pandit2020generalized, wang2022high, wang2023rate, chen2023community}, among others.

Although many mechanisms, such as regulatory processes in biology (cf. \cite{sima2009} for a survey), involve nonlinear dynamics, research on high-dimensional time series models addressing such dynamics remains limited. \cite{mazur2009} and \cite{aijo2009} employed Bayesian learning to manage the stochasticity of biological data. \cite{lim2015} introduced a family of VAR models using operator-valued kernels to identify nonlinear dynamic systems. \cite{zhou2018} proposed a framework for non-parametric autoregressive models within generalized linear models by utilizing reproducing kernel Hilbert spaces, analyzing the convex penalized sparse and smooth estimator. \cite{shen2019nonlinear} investigated nonlinear structural VAR models with application to brain networks. Additional applications can be found in \cite{pereda2005nonlinear, balcilar2016testing, yu2021sparse}, among others. Among these works, only \cite{zhou2018} provided theoretical guarantees, although their concentration inequalities are not sharp. In this paper, we extend the framework of sparse linear VAR models to sparse non-parametric nonlinear VAR models, with rigorous theoretical guarantees.

This paper has two primary objectives: (i) to develop sharp inequalities for tail probabilities for non-sub-Gaussian nonlinear VAR processes; (ii) to propose a new class of methods for high-dimensional non-parametric VAR models and to apply our inequalities to obtain theoretical properties of $\ell_1$ regularized estimators. It is expected that our framework, inequalities and tools will be useful in other high-dimensional linear and nonlinear VAR problems.

In our theoretical framework, we shall consider the following nonlinear VAR models
\begin{align} \label{eq:autoregmodel0}
X_t=h^{(1)}(X_{t-1})+h^{(2)}(X_{t-2})+\ldots+h^{(d)}(X_{t-d})+\epsilon_t,  
\end{align}
where $\epsilon_t\in\RR^p, t\in\ZZ$, are i.i.d. random vectors, $X_t=(X_{t}^{(1)},\ldots,X_{t}^{(p)})^\top\in\RR^p$, $h^{(j)}=(h_{1}^{(j)},\ldots,h_{p}^{(j)})^\top$ and $h_k^{(j)}:\RR^p\rightarrow \RR$, $1\le j\le d, 1\le k\le p,$ are real-valued functions. By stacking lagged vectors, we can let $d=1$ in \eqref{eq:autoregmodel0} and consider the nonlinear VAR(1) model. Then \eqref{eq:autoregmodel0} can be rewritten as
\begin{align} \label{eq:autoregmodel}
X_t=h(X_{t-1})+\epsilon_t.  
\end{align}
Based on model \eqref{eq:autoregmodel}, we shall develop sharp Bernstein-type inequalities. 
Establishing exponential-type tail probability inequalities for temporal dependent processes is a challenging problem. There has been some effort to derive concentration inequalities for non-i.i.d. processes. For example, generalizations of Bernstein's inequality to $\alpha$-mixing and $\phi$-mixing random variables have been studied in \cite{bosq1993}, \cite{modha1996}, \cite{samson2000} and \cite{merlevede2009, merlevede2011}, among others. \cite{zhang2021robust} provided Bernstein-type inequality for dependent random variables under geometric moment contraction. Exponential-type inequalities were also derived for sums of Markov chains in \cite{douc2008, adamczak2008, lemanczyk2021general}.
Unfortunately, all these inequalities involve extra non-constant factors to account for weak dependence, and are not as sharp as the original Bernstein's inequality for independent random variables. Recently, \cite{fan2021hoeffding} and \cite{jiang2018} established sharp Hoeffding-type inequality and Bernstein-type inequality for stationary Markov dependent random variables. \cite{chen2017} derived exponential inequalities and Nagaev-type inequalities for one dimensional linear (or moving average) processes under both short- and long-range dependence. 
Due to the interactions between temporal and cross-sectional dependence, tail probabilities of high-dimensional time series is much more complicated than the one-dimensional ones. 
In this work, we establish Bernstein-type inequalities for nonlinear VAR processes. Our inequalities, up to some constants, are as sharp as the classical Bernstein inequality for i.i.d. random variables. To the best of our knowledge, we are among the first to develop such sharp Bernstein-type inequalities for time series. Notably, we do not use the commonly employed ``blocking'' technique for sequences of dependent random variables \citep{hall2018}, which allows us to avoid logarithmic factors. Our technical approach can be used to improve existing studies on high-dimensional VAR models, such as in \cite{kock2015oracle, jiang2023autoregressive, dahlhaus2023adaptation, wang2023rate}.

To study nonlinear dynamical systems from high-dimensional time series data, in this paper, we introduce sparse additive non-parametric VAR models. Our method combines ideas from sparse linear modelling, additive non-parametric regression and VAR models. Each nonlinear function $h_j$, $1\le j\le p$, in model \eqref{eq:autoregmodel} can be expressed as: $$ h_{j}(x)=\sum_{k=1}^p h_{jk}(x_k),$$
where $x=(x_1,\ldots,x_p)^\top \in \RR^p$ and $h_{jk}(\cdot)$ are functions of one dimensional variables. The underlying VAR model is similar to sparse linear regression, but we impose a sparsity constraint on the index set $\{(j,k): h_{jk}(\cdot)\neq 0\}$ of functions $h_{jk}$ that are not identically zero. Then we estimate each nonlinear function $h_{jk}$ in terms of a truncated set of basis functions. \cite{ravikumar2009} proposed sparse additive linear models using a basis expansion and LASSO type penalty under i.i.d. data. \cite{meier2009} considered a sparsity-smoothness penalty for high-dimensional generalized additive models. \cite{koltchinskii2010}, \cite{raskutti2012minimax} and \cite{yuan2016} studied a different framework, sparse additive kernel regression, for the cases where the component functions belong to a reproducing kernel Hilbert spaces (RKHS). They penalized the sum of the reproducing kernel Hilbert space norms of the component functions. Their sparse additive linear models are extended to autoregressive generalized linear models in \cite{zhou2018}. \cite{lim2015} introduced operator-valued kernel-based VAR models, and developed proximal gradient descent algorithms. However, their paper does not provide any theoretical guarantees. Recently, \cite{duker2025kernel} developed an RKHS-based framework for nonlinear VAR processes and derived non-asymptotic probabilistic bounds.

In this work, our method has the nice feature that it decouples smoothness and sparsity. This leads to a simple block coordinate descent algorithm (cf. \cite{ravikumar2009}) that can be carried out with any non-parametric smoother and scales easily to high-dimensions. Besides, with our new probability inequalities as primary tools, we can analyze the properties of $\ell_1$ regularized estimators under non-Gaussian errors in the context where $p$ is much larger than $n$. Roughly speaking, $p$ can be as large as $e^{n^c}$ for some constant $0<c<1$ if $\epsilon_t$ has finite exponential moments, and the power constant $c$ is related to the truncated number of basis expansion. We shall give a detailed description on how the dependence measures of the processes, the moment condition of the errors, the sparsity of functions and basis expansion affect the rate of convergence and the model selection consistency of the estimator. 

The rest of the paper is structured as follows. Section \ref{section:bern} presents Bernstein-type inequalities for nonlinear VAR processes in \eqref{eq:autoregmodel} under Lipschitz condition and different types of moment conditions for the error processes. In Section \ref{section:var}, we first formulate an $\ell_1$ regularized optimization problem for nonlinear VAR models on the population level that induces sparsity. Then we derive a sample version of the problem using basis expansion. Theoretical properties that analyze the effectiveness of the estimators in the high-dimensional setting are also presented. Simulation studies and real data analysis are carried out in Sections \ref{section:simulation} and \ref{section:data}, respectively. Proofs of theorems and technical lemmas are contained in Section \ref{section:proofs}.

We now introduce some notation. For a vector $x=(x_1,\ldots,x_p)^\top$, define $\|x\|_q = (|x_1|^q+\ldots+|x_p|^q)^{1/q}$, $q\ge 1$, $\|x\|=\|x\|_2$, $\|x\|_{\infty}=\max_{1\le j\le p} |x_j|$, and 
$\abs(x):=(|x_1|,\ldots,|x_p|)^\top.$ For a matrix $A = (a_{ij})$, write $|A|_{\infty}=\max_{i,j}|a_{ij}|$, the Frobenius norm $\|A\|_F = (\sum_{ij} a_{ij}^2)^{1/2}$, the spectral norm $\|A\|_2=\max_{\|x\|_2\le 1} \|Ax\|_2$ and the matrix infinity norm $\|A\|_{\infty}=\max_{i}\sum_j |a_{ij}|$. Let $\lambda_{\min}(A)$ (resp. $\lambda_{\max}(A)$) be the minimum (resp. maximum) eigenvalue of $A$. 
For two sequences of real numbers $\{a_n\}$ and $\{b_n\}$, write $a_n=O(b_n)$ (resp. $a_n\asymp b_n$) if there exists a constant $C$ such that $|a_n|\leq C |b_n|$ (resp. $1/C \leq a_n/b_n\leq C$) holds for all sufficiently large $n$, and write $a_n=o(b_n)$ if $\lim_{n\to\infty} a_n/b_n =0$.

Let $\epsilon_t, t\in\mathbb Z$, be i.i.d. random vectors and $\F_k=(\ldots,\epsilon_{k-1},\epsilon_k)$. Define projection operator $P_k$, $k\in\ZZ$, by $P_k(\cdot)=\EE(\cdot|\F_k)-\EE(\cdot|\F_{k-1}).$ Let $(\epsilon_k')$ be an i.i.d. copy of $(\epsilon_k)$. For $X_t=g(\ldots,\epsilon_{t-1},\epsilon_t)$, where $g$ is a measurable function, we define the coupled version $X_{t,\{k\}}=g(\ldots,\epsilon_{k-1},\epsilon_{k}',\epsilon_{k+1},\ldots,\epsilon_t)$, which has the same distribution as $X_t$ with $\epsilon_{k}$ in the latter replaced by an i.i.d. copy $\epsilon_{k}'$.

\section{Bernstein Inequalities for Nonlinear VAR Processes}\label{section:bern}

Exponential inequalities play a fundamental role in high-dimensional inference. Differently from i.i.d. random variables, directly applying concentration inequalities for dependent random variables to high-dimensional time series problems may lead to suboptimal results in many cases, due to the interrelationship between temporal and cross-sectional dependencies. \cite{zhang2017,zhang2021convergence,han2023high} introduced new dependence measures to describe temporal and cross-sectional dependence of high-dimensional time series, then derived Fuk-Nagaev type inequalities for heavy tailed random vectors to study statistical properties of sample mean vector, spectral density matrix estimation and robust $M$-estimation, respectively. In this section, we shall present new and powerful inequalities for tail probabilities of nonlinear vector autoregressive (VAR) processes. The processes can be non-Gaussian. In Theorem \ref{thm:bern}, we provide Bernstein-type inequalities for nonlinear VAR process under finite moment condition and exponential moment condition, respectively. In contrast, exponential inequalities provided in \cite{basu2015} are only applicable to Gaussian processes and linear VAR models with Gaussian innovation vectors (cf. Proposition 2.4 therein).

To establish exponential inequalities, we introduce the following assumptions on the function $h$ and the errors $\epsilon_t$ in model \eqref{eq:autoregmodel}. Recall that $\| \cdot\|_\infty$ is the matrix infinity norm.

\begin{assumption}\label{asmp:LipconstH}
Consider model \eqref{eq:autoregmodel}, let $h=(h_1,\ldots,h_p)^\top$ and $h_j: \mathbb R^p \rightarrow \mathbb R$, $1\le j\le p$ be real valued functions. Assume that componentwise Lipschitz condition holds for each $h_j$. That is, for any $x=(x_1,\ldots,x_p)^\top,y=(y_1,\ldots,y_p)^\top\in\mathbb R^p$, $1\le j\le p$, there exist coefficients $H_{jk}\ge 0$ such that
\begin{equation} \label{eq:defH}
|h_j(x) -h_j(y)| \le \sum_{k=1}^p H_{jk} |x_k-y_k|.
\end{equation}
Write $H=(H_{jk})_{p \times p}$ and $\|H\|_\infty=\max_{1\le j\le p}\sum_{k=1}^p H_{jk}$. 
Assume there exists an absolute constant $0<\rho<1$ such that $\|H\|_\infty\leq \rho$. 
\end{assumption}

The above assumption requires componentwise Lipschitz condition for nonlinear VAR processes. This assumption can be easily extended to nonlinear VAR($d$) processes. See also \cite{chen1993}, \cite{diaconis1999}, \cite{jarner2001}, \cite{shao2007}, \cite{fan2008} and \cite{chen2016} for nonlinear autoregressive processes. 
Intuitively, $\rho$ quantifies the strength of dependence. For example, in one dimensional AR(1) model, $X_t=\rho X_{t-1}+\epsilon_t$. Larger $\rho$ suggests stronger dependence.

\begin{remark}[Existence of stationary distribution] For the sake of completeness, in this remark, we shall apply the theory in \cite{chen2016} and show the existence of stationary distribution.
Construct a collection of backward series of random vectors $X_{(-n),t}$, for $t \geq -n$, as follows. For all $t \in \mathbb{Z}$, define
$X_{(t),t} = 0$ and the recursion,
\[
X_{(-n), t} = h\big(X_{(-n), t-1}\big) + \epsilon_t.
\]
Let $X_{(-n),t}^{(j)}$ denote the $j$-th component of $X_{(-n), t}$. Then $X_{(-n), t}
=(X_{(-n),t}^{(1)},\ldots,X_{(-n), t}^{(p)})^\top$. 
Under Assumption \ref{asmp:LipconstH}, we have
\begin{align}
\label{eq:stationaryxtjn}
\|X_{(-n+1), t}-X_{(-n), t}\|_\infty  
&\leq \max_{1\leq j\leq p}\sum_{k=1}^p H_{jk}|X_{(-n+1), t-1}^{(k)}-X_{(-n), t-1}^{(k)}| \nonumber\\
&\leq \rho \|X_{(-n+1), t-1}-X_{(-n), t-1}\|_\infty\nonumber\\
&\leq \rho^{t+n-1}
\|X_{(-n+1),-n+1}-X_{(-n),-n+1}\|_\infty. 
\end{align}
Taking the $L_q$ norm and defining $c_p=\|\max_{1\leq j\leq p}|h_j(0)+\epsilon_{1}^{(j)}|\|_q$, we obtain
\begin{align*}
\big\|\|X_{(-n+1),t}-X_{(-n),t}\|_\infty\big\|_q\leq \rho^{t+n-1}
c_p.    
\end{align*}
Since $\rho<1$, for fixed $t$, the sequence $X_{(-n), t}^{(j)}$ converges as $n\rightarrow\infty$ for any $1\leq j\leq p$. Denote the limit by $Y_{t}^{(j)}$ and set $Y_t=(Y_{t}^{(1)},\ldots,Y_{t}^{(p)})^\top.$ We now show that $Y_t$ is the stationary solution of model \eqref{eq:autoregmodel}. For any $\kappa>0,$ there exists an $N_0\in\NN$, such that for $n>N_0$ we have $\|Y_{k}-X_{(-n),k}\|_\infty<\kappa$. Then, using the Lipschitz condition on $h$,
\begin{align*}
\max_{1\leq j\leq p}|h_j(Y_{t-1})+\epsilon_{t}^{(j)}-Y_{t}^{(j)}|
\leq 
\max_{1\leq j\leq p}|h_j(Y_{t-1})-h_j(X_{(-n),t-1})|
+\max_{1\leq j\leq p}|X_{(-n),t}^{(j)}-Y_{t}^{(j)}| \leq 2\kappa. 
\end{align*}
Since $\kappa$ is arbitrary, this implies $Y_t=h(Y_{t-1})+\epsilon_t$ almost surely, so $Y_t$ indeed satisfies the VAR recursion and is stationary.

\end{remark}

\begin{remark}
Denote $H^m$ as the result of multiplying the matrix $H$ by itself $m$ times. We can weaken Assumption~\ref{asmp:LipconstH} by requiring only that there exists an integer
\( m \ge 1 \) such that \( \|H^m\|_\infty \leq \rho < 1 \). In other words, we allow \( \|H\|_\infty \geq 1 \) as long as repeated application of \( H \) eventually satisfies this condition. 
All of our subsequent results remain valid under this relaxed assumption.
To see this, we consider the previous argument for the existence of a stationary distribution. Repeatedly applying the first inequality in \eqref{eq:stationaryxtjn}, we have
 \begin{align}
 \label{eq:xtjtn1}
\|X_{(-n+1),t}-X_{(-n),t}\|_\infty  
&\leq \max_{1\leq j\leq p}\sum_{k=1}^p (H^m)_{jk}|X_{(-n+1),t-m}^{(k)}-X_{(-n),t-m}^{(k)}| \nonumber\\
&\leq \rho \|X_{(-n+1),t-m}-X_{(-n),t-m}\|_\infty.
\end{align}
Iterating \eqref{eq:xtjtn1}, we conclude that 
$\|X_{(-n+1),t}-X_{(-n),t}\|_\infty\lesssim \rho^{\lfloor (n+t)/m\rfloor},$ where $\lfloor x\rfloor$ is the largest integer less or equal to $x.$
Then for fixed $t$, we have that $X_{(-n),t}$ converges as $n\rightarrow \infty$. Similar adaptations apply throughout the paper under this relaxed condition instead of Assumption \ref{asmp:LipconstH}.
\end{remark}

\begin{assumption}
\label{asmp:moment}
For i.i.d. random vectors $\epsilon_t = (\epsilon_{t}^{(1)}, \ldots, \epsilon_{t}^{(p)})^\top \in\RR^p,$ $t\in\ZZ,$ assume one of the following holds:
\begin{enumerate}
    \item [(i)](finite moment) $\mu_q:=\max_{1\leq j\leq p} (\EE |\epsilon_{t}^{(j)}|^q )^{1/q}<\infty$ for some $q\geq 2.$ 
    \item[(ii)](exponential moment) $\mu_e:=\max_{1\le j\le p}\EE\big(\exp(c_0 |\epsilon_{t}^{(j)}| )\big)$, for some $c_0>0$. 
\end{enumerate}
\end{assumption}

\begin{assumption}
\label{asmp:functiong}
Let function $g:\RR^p\rightarrow \RR,$ be Lipschitz continuous with $|g(x)-g(y)|\leq \sum_{j=1}^p G_j |x_j-y_j|,$ for any $x=(x_1,\ldots,x_p)^\top,$ $y=(y_1,\ldots,y_p)^\top\in\RR^p$, where $G_j$ are Lipschitz coefficients. Denote $G=(G_1,\ldots,G_p)^\top$ and $\tau:=\|G\|_1=\sum_{j=1}^p G_j$. 
\end{assumption}

The following theorem presents a Bernstein-type inequality for bounded Lipschitz continuous functions, under both the finite moment condition and the exponential moment condition of the error vectors $\epsilon_t$, respectively.


\begin{theorem} \label{thm:bern}
Consider the VAR process defined in \eqref{eq:autoregmodel}, where the function \( h \) satisfies Assumption~\ref{asmp:LipconstH}. Let \( g \) be any function satisfying Assumption~\ref{asmp:functiong} with $\tau=\|G\|_1$. Then:
\begin{itemize}
    \item[(i)] If Assumption~\ref{asmp:moment} $(i)$ holds and \( g \) is bounded with \( \|g\|_\infty =\sup_x |g(x)|\le M \), then for all \( z \ge 0 \),
    \begin{align} \label{inequality:bern}
    \PP\left(\left|\sum_{t=1}^n \left(g(X_t) - \EE g(X_t)\right)\right| \ge z\right)
    \le 2 \exp\left\{-\frac{z^2}{c_1 \tau^2 n + c_2 \tau M z}\right\},
    \end{align}
    where \( c_1 \) and \( c_2 \) are positive constants depending only on \( q \), \( \rho \), and \( \mu_q \).    \item[(ii)] If Assumption~\ref{asmp:moment} $(ii)$ holds, then for all \( z \ge 0 \),
    \begin{align} \label{inequality:bern2}
    \PP\left(\left|\sum_{t=1}^n \left(g(X_t) - \EE g(X_t)\right)\right| \ge z\right)
    \le 2 \exp\left\{-\frac{z^2}{c_3 \tau^2 n + c_4 \tau z}\right\},
    \end{align}
    where \( c_3 \) and \( c_4 \) are positive constants depending only on \( \rho \) and \( \mu_e \).
\end{itemize}
\end{theorem}

Theorem \ref{thm:bern}(i) addresses the finite moment case for the error vectors $\epsilon_t$ (cf. Assumption \ref{asmp:moment} (i)). 
If the error vectors $\epsilon_t, t\in\ZZ,$ satisfy stronger moment condition than merely having a finite $q$-th moment, we can expect a stronger inequality than \eqref{inequality:bern}. Indeed, when $\epsilon_t$ has subexponenial tail (Assumption \ref{asmp:moment} (ii)), we obtain an improved Bernstein-type inequality in \eqref{inequality:bern2}. Different from Theorem \ref{thm:bern} (i), in Theorem \ref{thm:bern}(ii), function $g$ can be unbounded.

\begin{remark}
Based on the proof of Theorem \ref{thm:bern}(i), we can have the explicit form for coefficients $c_1$ and $c_2$ as $c_1=32e^2(-\rho^2\log\rho)^{-2}\mu_2^2$ and $c_2=8e(-\rho^2\log\rho)^{-1}$. If function $g$ is bounded by an absolute constant, then we can simplify above tail inequality \eqref{inequality:bern} and obtain the following Hoeffding type inequality.
\end{remark}
\begin{corollary}
Consider the VAR process defined in \eqref{eq:autoregmodel}, where the function \( h \) satisfies Assumption~\ref{asmp:LipconstH}. Let \( g \) be any function satisfying Assumption~\ref{asmp:functiong}. Suppose Assumption~\ref{asmp:moment} $(i)$ or \ref{asmp:moment} $(ii)$ holds. 
If $g$ is bounded with $\|g\|_\infty\leq 1$, then we have
\begin{align} \label{inequality:hoeff}
\PP\Big(\Big|\sum_{t=1}^n\big(g(X_t)-\EE g(X_t)\big)\Big|\geq z\Big)\leq 2e^{-c_1 z^2/(\tau^2n)}, 
\end{align}
where $c_1$ is a positive constant depending only on $q$, $\rho$ and $\mu_q$.
\end{corollary}

\begin{remark}
Note that up to a multiplicative constant, our Bernstein-type inequality \eqref{inequality:bern} coincides with classical Bernstein's inequality for i.i.d. random variables. Thus one can expect sharper convergence rates for estimators of nonlinear VAR processes \eqref{eq:autoregmodel}. We remark that the majority of the previous inequalities for temporal dependent processes do not recover Bernstein's inequality. For example, under geometric moment contraction with decay coefficient $0<\rho<1$ (see \cite{wu2004}) and assume $|X_t|\le M$, \cite{zhang2021robust} provided the following Bernstein-type inequality,
\begin{align*}
\PP\Big(\left|\sum_{t=1}^n \big( X_t-\EE X_t\big)\right|\geq z\Big)\leq \exp\left\{-\frac{z^2}{4c_1(c_3n+M^2)+2c_2M(\log(n))^2 z} \right\}, 
\end{align*}
where $c_1, c_2$ are some constants only depending on $\rho$, and $c_3<\infty$ is a positive constant measuring the temporal dependence. Similarly, \cite{merlevede2009} obtained a Bernstein-type inequality for a class
of exponentially decay $\alpha$-mixing and bounded random variables,
\begin{align*}
\PP\Big(\left|\sum_{t=1}^n \big( X_t-\EE X_t\big)\right|\geq z\Big)\leq \exp\left\{-\frac{c_1 z^2}{nM^2+M\log(n)\log\log(n) z} \right\}, 
\end{align*}
where $c_1>0$ and $|X_t|\le M$. Both involve an unpleasant $\log(n)$-type multiplicative factor. Our sharp Bernstein-type inequality is of independent interest. 
We expect our sharp inequality can be useful for other high-dimensional linear and nonlinear time series problems. 
\end{remark}

\paragraph{Proof Sketch.}
The proof of Theorem \ref{thm:bern} is quite involved. The key steps involve employing a martingale decomposition and deriving a sharp bound for the martingale differences. To be more specific, without loss of generality, assume $\|G\|_1=1$ with $G$ defined in Assumption~\ref{asmp:functiong}. Recall that $\F_k=(\ldots,\epsilon_{k-1},\epsilon_k)$ and the projection operator $P_k(\cdot)=\EE(\cdot|\F_k)-\EE(\cdot|\F_{k-1})$, for $k\in\ZZ$.  
The summation can then be decomposed into a sum of martingale differences:
\begin{align*}
S_n(g):=\sum_{t=1}^n \big(g(X_t)-\EE g(X_t)\big)
=\sum_{k\leq n}\xi_k, \quad\textrm{where } 
\xi_k=P_k(S_n(g)).
\end{align*}
For $X_t=\mathcal G(\ldots,\epsilon_{t-1},\epsilon_t)$, where $\mathcal G$ is some measurable function, following \cite{wu2005nonlinear}, we define the coupled version $$X_{t,\{k\}}=\mathcal G(\ldots,\epsilon_{k-1},\epsilon_{k}',\epsilon_{k+1},\ldots,\epsilon_t).$$ For $x=(x_1,\ldots,x_p)^\top$, write $\abs(x)=(|x_1|,\ldots,|x_p|)^\top$.
Since the mapping $h$ is componentwise Lipschitz continuous, by induction, we have $\abs(X_t-X_{t,\{k\}})\leq H^{t-k} \abs(\epsilon_k-\epsilon_k')$. Hence 
\begin{align}
\label{eq:bddofpkgxt}
\big|P_k (g(X_t))\big|
&=\big|\EE(g(X_t)-g(X_{t,\{k\}})|\F_k)\big|\nonumber\\
&\leq \EE \left( G^\top \abs(X_t-X_{t,\{k\}})  \big|\F_k\right)\nonumber\\
&\leq  \EE \left( G^\top  H^{t-k} \abs(\epsilon_k-\epsilon_k') \big|\F_k\right).
\end{align}
Since the function \( g(\cdot) \) is bounded by \( M \), it follows that \( |P_k(g(X_t))| \leq 2M \). Therefore, combining this with \eqref{eq:bddofpkgxt}, we obtain
\begin{align}
|\xi_k|\leq \sum_{t=1}^n |P_k(g(X_t))|\leq \sum_{t=k\vee 1}^n\min\Big\{v_{t-k}^\top\EE(\abs(\epsilon_k-\epsilon_k') \big |\F_k), 2M \Big\},\quad\mathrm{with}\quad \|v_t\|_1\leq \rho^t.  
\end{align}
Since $\|v_t\|_1$ decays exponentially fast, for all sufficiently large $t$, one shall expect the first term $v_{t-k}^\top\EE(\abs(\epsilon_k-\epsilon_k') \big |\F_k)$ to be small.
Then by carefully leveraging between the two terms as detailed in Lemma \ref{lem:increaseorder}, we obtain that $$\EE(e^{|\xi_k|h})<\infty$$ for any $h\leq h^*$ some constant $h^*>0$.   
Since $\xi_k$'s are martingale differences, 
\begin{align}
\EE(e^{\xi_kh}|\F_{k-1})
&=1+\EE(e^{\xi_kh}-\xi_kh-1|\F_{k-1})\nonumber\\
&\leq 1+\EE\Big[\frac{e^{|\xi_k|h}-|\xi_k|h-1}{h^2}\Big|\F_{k-1}\Big]h^2,
\end{align}
where the conditional expectation in the last line can be shown to be bounded for any $h\leq h^*$ with the bound denoted by $c$. Hence 
\begin{align}
\label{eq:bddexp}
\EE(e^{\xi_kh}|\F_{k-1})\leq 1+c h^2.
\end{align}
The above applies for $-n\leq k\leq n.$ For $k<-n$, we can show that those terms are negligible.
The desired result then follows by Markov's inequality and recursively using \eqref{eq:bddexp} for $-n\leq k\leq n.$

It should be emphasized that our Bernstein-type concentration inequalities are sharp, and does not contain any annoying extra logarithmic terms. These inequalities are useful for handling non-Gaussian VAR problems. 

\section{Sparse additive nonlinear VAR models}\label{section:var}

In this section, we study sparse additive nonlinear VAR models. We first introduce a population-level optimization problem and then derive a sample-based algorithm through basis expansion. Our theoretical analysis builds on the technical tools developed in previous section.

\subsection{The model} \label{section:model}

Assume that we are provided with observed time series data $X_1,\ldots,X_n \in \RR^p$, which are sampled from a dynamical system involving $p$ variables. Our primary goal is to infer the direct influence that each variable $j$ exerts on every other variable $k$ (with $k\neq j$, $1\le k \le p$). For instance, in the case of linear VAR models, the evolution of the system is often characterized by $X_t=GX_{t-1}+\epsilon_t$, where $G$ is a $p\times p$ coefficient matrix, and $\epsilon_t$ represents noise. In our study, we assume that a first-order stationary model provides a sufficient approximation of the temporal dependencies within the system. Accordingly, we recall the nonlinear VAR framework in \eqref{eq:autoregmodel},
\begin{align*}
X_t=h(X_{t-1})+\epsilon_t,  
\end{align*}
where the function $h$ can capture potentially complex, nonlinear dynamics.

In this section, we propose a new class of high-dimensional, sparse, additive non-parametric VAR models. Here, each component $h_j$ of the function $h$ is assumed to decompose additively in terms of the individual components of the state vector $x\in\RR^p$. Specifically, we posit that for each variable $j$
\begin{align}
\label{eq:additivefunc}
h_j(x)=\sum_{k=1}^p h_{jk}(x_k),
\end{align}
where each function $h_{jk}:\RR\rightarrow \RR$ captures the individual contribution of the $k$-th variable to the dynamics of the $j$-th variable.


Let $\Pi$ denote the joint distribution of the vector $X_t$, and let $\Pi_k$ denote the marginal distribution of the $k$-th component $X_{t}^{(k)}$ for each $1\le k\le p$. For practical purposes, we define the $L_2(\Pi_k)$-norm of the function $h_{jk}$ as
\begin{eqnarray*}
\| h_{jk} \|_{\Pi_k,2} =\sqrt{ \int h_{jk}^2(x) \mathrm{d} \Pi_k(x) } = \sqrt{\EE h_{jk}^2(X_{t}^{(k)})}.
\end{eqnarray*} 
This definition is particularly relevant because it allows us to accommodate functions $h_{jk}$ that might not be Lebesgue integrable over the entire real line; instead, the integrability is considered relative to the distribution $\Pi_k$.


The classical nonlinear ridge regression is defined as 
\begin{align*} 
\frac1n\sum_{t=1}^n \| X_t-h(X_{t-1})  \|_2^2 + \lambda\sum_{j=1}^p \sum_{k=1}^p \| h_{jk}\|_{\Pi_k,2}^2 ,
\end{align*}
where the norms measure the overall discrepancy and the smoothness penalty on each component.
To encourage sparsity in high-dimensional settings, we replace the squared norm $\| h_{jk}\|_{\Pi_k,2}^2$ with the norm $\| h_{jk}\|_{\Pi_k,2}$ itself. This substitution leads to a population-level penalized least squares estimator defined by the optimization problem
\begin{eqnarray} \label{problem}
(\hat h_{jk},1\leq j,k\leq p):= \underset{h_{jk}\in\mathcal I_{k},1\le j,k\le p }{\argmin} \left\{\frac1n\sum_{t=1}^n \| X_t-h(X_{t-1})  \|_2^2 + \lambda\sum_{j=1}^p\sum_{k=1}^p \| h_{jk}\|_{\Pi_k,2}\right\} .
\end{eqnarray}
Here, $h$ is decomposed as in \eqref{eq:additivefunc} into a sum of univariate functions, and $\mathcal I_{k}$ is an appropriate function class for the $k$-th component. In practice, the norm $\| h_{jk}\|_{\Pi_k,2}$ can be estimated empirically by $(n^{-1}\sum_{t=1}^nh_{jk}^2(X_{t-1}^{(k)}))^{1/2}$.

By decomposing $h_{j}$ into additive components, our framework enhances interpretability and computational efficiency in high-dimensional settings. The imposed sparsity helps to pinpoint which variables have a direct influence on the dynamics. This model is an extension of the sparse additive models developed for the i.i.d. case \citep{ravikumar2009} and is especially relevant when the system exhibits nonlinear structure that traditional linear models fail to capture, while still preserving a structure that is amenable to rigorous analysis and estimation.

For each $k\in\{1,\ldots,p\}$, let $\mathcal H_k$ denote the Hilbert subspace $L_2(\Pi_k)$ consisting of measurable functions $f(\cdot)$ satisfying $\EE f(X_t^{(k)})=0$ and the norm $\|f\|= (\EE f^2(X_t^{(k)}) )^{1/2} < \infty$. The inner product on $\mathcal H_k$ is defined as 
$$\langle f,g \rangle=\EE\big( f(X_t^{(k)}) g(X_t^{(k)})\big).$$ 
We denote by $\mathcal H=\mathcal H_1\oplus \mathcal H_2\oplus \cdots \oplus H_p$ the Hilbert space of functions of $(x_1,\ldots,x_p)$ that admit an additive representation $m(x)=\sum_k f_k(x_k)$ with each $f_k\in\mathcal H_k, k=1,\ldots,p$.

We now impose the following assumption on our basis expansion.

\begin{assumption}[Basis function] \label{asmp:expansion}
Assume that the functions $h_{jk}(x)$ in \eqref{eq:additivefunc} have compact support for all $1 \leq j, k \leq p$, that is, 
$|h_{jk}(x)| = 0$ for any $|x| > c_0$, for some constant $c_0 > 0$. Moreover assume $h_{jk}\in \mathcal I_{k}$ where
\begin{align*}
\mathcal I_{k} = \left\{ h_{jk}(\cdot)\in\mathcal H_k: h_{jk}(\cdot)=\sum_{l=1}^\infty b_{jk}^{(l)*} \psi_{k,l}(\cdot),\quad \sum_{l=1}^\infty (b_{jk}^{(l)*})^2 l^{2\beta} \le C^2\right\},
\end{align*}  
where $(\psi_{k,l}(\cdot): l=1,2,\ldots)$ is a uniformly bounded orthonormal basis on $[-c_0,c_0]$, that is $|\psi_{k,l}(x)|\leq B,$
for some $0<B,C<\infty$ and $\beta\ge 1$. 
\end{assumption}

For example, we can choose the Fourier basis functions to satisfy Assumption \ref{asmp:expansion}. In standard nonparametric regression such as \cite{ravikumar2009}, covariates are often assumed to be bounded (i.e., to have compact support). Similarly, in our nonlinear VAR framework we assume $h_{jk}(x)$ in \eqref{eq:additivefunc} have compact support for mathematical convenience and tractability; see also \cite{raskutti2012minimax,zhou2018}. Many of our results can be extended to the case of unbounded $h_{jk}(x)$ via truncation arguments with proper tail decay conditions. We omit such arguments for the sake of presentation simplicity.
For example, the Fourier basis satisfies this assumption.
This assumption implies that the tail of the expansion satisfies $\sum_{l=L+1}^\infty (b_{jl}^{(l)*})^2 \le C^2L^{-2\beta}$, which corresponds to the functional class condition of \cite{ravikumar2009} and is a standard requirement in basis expansion methods. The parameter $\beta$ captures the level of smoothness, effectively linking our function class to a function space. Although one could allow $\beta$ to vary adaptively with $k$, we confine ourselves to a common smoothness level in this work.

Let $L=L_n$ be a truncation parameter, and let $h_{jk}^{(L)}$ be the approximation of $h_{jk}$ defined by 
\begin{align}
h_{jk}^{(L)}(\cdot)=\sum_{l=1}^L b_{jk}^{(l)*}\psi_{k,l}(\cdot).    
\end{align}
In this formulation, $h_{jk}^{(L)}$ is interpreted as the projection of $h_{jk}$ onto the truncated set of basis functions $\{\psi_{k,1}, \ldots, \psi_{k,L}\}$.
Then, for $1\le j,k\le p,$ the model can be written as
\begin{eqnarray}
\label{eq:xij1}
X_{t}^{(j)} = \sum_{k=1}^p h_{jk}^{(L)}(X_{t-1}^{(k)})+r_{t}^{(j)}+\epsilon_{t}^{(j)}, \mbox{ where }
 r_{t}^{(j)}=\sum_{k=1}^p [h_{jk}(X_{t-1}^{(k)})- h_{jk}^{(L)}(X_{t-1}^{(k)})]
\end{eqnarray}
is the reminder term and captures the bias introduced by truncating the basis expansion.

We now define the oracle coefficients on the population level for the basis expansion and the design matrix. For any $x=(x_1,\ldots,x_p)^\top \in \RR^p$, set vectors
\begin{equation} \label{eq:def:b}
\begin{split}
& b_{j,k}^*=(b_{j,k}^{(1)*},\ldots,b_{j,k}^{(L)*})^\top, \\
& b_{j}^*=(b_{j,1}^{*\top},\ldots,b_{j,p}^{*\top})^\top, \\
& b^*=(b_{1}^{*\top},\ldots,b_{p}^{*\top})^\top, \\
& \psi_{k} (x_k) = (\psi_{k,1}(x_k),\ldots,\psi_{k,l}(x_k))^\top, \\
& \psi (x) = (\psi_{1}^\top (x_1) ,\ldots,\psi_{p}^\top (x_p) )^\top.
\end{split}
\end{equation}
Let $r_t=(r_{t}^{(1)},\ldots,r_{t}^{(p)})^\top$. With these definitions, the model can be rewritten in a compact form as 
\begin{align}
\label{eq:XiPhi}
 X_t&:= \Psi(X_{t-1})^\top b^* + r_t + \epsilon_t,
\end{align}
where
\begin{align*}
\Psi(X_{t-1})&=\left(\begin{array}{cccccc}
\psi (X_{t-1}) &  0 & 0 &\cdots  & 0 \\
0 & \psi (X_{t-1})  & 0 &\cdots  & 0\\
0 & 0 & \psi (X_{t-1})  &\cdots  & 0\\
\vdots & \vdots & \vdots & \ddots & \vdots\\
0 & 0 & 0 &\cdots  & \psi (X_{t-1}) \\
\end{array}\right) \in\RR^{p \times p^2 L}.
\end{align*}

Consequently, the solution to our optimization problem \eqref{problem} can be approximately estimated by solving
\begin{eqnarray} \label{problem2}
\hat b:= \underset{b }{\argmin} \left\{\frac1n \sum_{t=1}^n \| X_t-\Psi(X_{t-1})^\top b \|_2^2 + \lambda\sum_{j=1}^p\sum_{k=1}^p \sqrt{\frac1n \sum_{t=1}^n \left(\sum_{l=1}^L \psi_{k,l} (X_{t-1}^{(k)}) b_{j,k}^{(l)} \right)^2 }\right\}.
\end{eqnarray}
This formulation can be viewed as a functional version of the group lasso, and the standard convexity arguments guarantee the existence of a minimizer.

Compared with the approach in \cite{lim2015}, which employs operator-valued reproducing kernels for VAR models, our formulation offers a key advantage: it decouples the smoothness and sparsity components. This separation allows us to employ a block coordinate descent algorithm (cf. \cite{ravikumar2009}) to efficiently construct the estimator. In the following section, we leverage the technical tools developed in Section \ref{section:bern} to establish the theoretical properties of our $\ell_1$-regularized estimator, under the assumption that the particular smoother in \eqref{problem2} is used.

\subsection{Asymptotic properties}
To facilitate the theoretical analysis, we impose the following assumptions on the functions $h_{jk}$ ($1\le j,k\le p$) and the basis expansions. For a function $f:\RR^d\rightarrow\RR,$ denote $\|f\|_2:=(\int_{\RR^d} f^2(x)d x)^{1/2}$. 




\begin{assumption} \label{asmp:eig}
There exist constants $\phi_U, \phi_L>0$, so that 
\begin{align} \label{eq:eig1}
\lambda_{\min}\Big\{\EE \psi(X_{t-1})\psi(X_{t-1})^\top\Big\}\ge \phi_L,
\end{align}
and
\begin{align} \label{eq:eig2}
\max_{1\leq k\leq p}\lambda_{\max}\Big\{ \EE \psi_{k} (X_{t-1}^{(k)})\psi_{k} (X_{t-1}^{(k)})^\top\Big\} \le \phi_U.
\end{align}
\end{assumption}


Condition \eqref{eq:eig1} in Assumption \ref{asmp:eig} is similar to the smallest population eigenvalue conditions commonly used in high-dimensional statistics \citep{raskutti2011minimax,vandegeer2014}. In addition, it parallels the population minimum eigenvalue condition in Assumption 4 of \cite{chen2015optimal} and Assumption S.3 of \cite{belloni2019conditional} for sieve basis expansion functions.
If the marginal density of $X_t^{(k)}$ satisfies $0< f_{\min} \le f_k(x)\le f_{\max} <\infty$ for $1\le k\le p $ and almost all $x\in[-c_0,c_0]$, then
\begin{align*}
&\EE u^\top \psi_{k} (X_{t-1}^{(k)})\psi_{k} (X_{t-1}^{(k)})^\top  u  =  \int_{-c_0}^{c_0} (u^\top \psi_{k} (x) )^2 f_k(x)\; {\rm d} x \le f_{\max}\; \|u\|_2^2 , \\
&\EE u^\top \psi_{k} (X_{t-1}^{(k)})\psi_{k} (X_{t-1}^{(k)})^\top  u  =  \int_{-c_0}^{c_0} (u^\top \psi_{k} (x) )^2 f_k(x)\; {\rm d} x \ge f_{\min}\; \|u\|_2^2 .
\end{align*}
This verifies condition \eqref{eq:eig2} in Assumption \ref{asmp:eig}.
In the following Proposition \ref{prop:eig}, we use concentration inequalities to establish the sample version of Assumption \ref{asmp:eig}.

\begin{proposition}\label{prop:eig}
Suppose Assumptions \ref{asmp:LipconstH} and \ref{asmp:moment}(ii) hold. Assume $\sup_x|\psi_{k,l}(x)|\le B$ for any $1\le k\le p, 1\le l\le L$. 

(i). Assume that \eqref{eq:eig1} holds
and that for some constant $c_1 >0$ does not rely on $p,L$, such that for all $u\in\mathbb R^{p L}$, 
\begin{align}\label{eq:thml2:moment}
\EE(u^\top\psi(X_t)\psi(X_t)^\top u)^2\le c_1 \big(u^\top\EE(\psi(X_t)\psi(X_t)^\top )u\big)^2.    
\end{align}
Then, with probability at least $1-p^{-c_2}-p e^{-c_3 n/\log(n)}$, for all $u\in\mathbb R^{pL}$ with $\|u\|_2=1$,
\begin{align}\label{eq:eiga}
\frac1n \sum_{t=1}^n u^\top \psi (X_t) \psi (X_t)^\top u \ge \frac{\phi_{L}}{2}-\frac1n  - c_4 \frac{\log(n)\log (pL)\cdot \|u\|_1^2}{n},
\end{align}
where $c_2,c_3,c_4>0$ are constants independent of $n,p,L$.

(ii). Assume that \eqref{eq:eig2} holds. 
Then, with probability at least $1-p^{-c_5}-e^{-c_6 n/\log(n)}$, for all $u\in\mathbb R^{L}$ with $\|u\|_2=1$,
\begin{align} \label{eq:eigb}
\max_{1\leq k\leq p}  \frac{1}{n}\sum_{t=1}^n u^\top  \psi_{k} (X_{t-1}^{(k)})\psi_{k} (X_{t-1}^{(k)})^\top u \le \phi_{U}+ c_{7} L\sqrt{\frac{\log(n)(\log p+\log L)}{n}},
\end{align}
where $c_5,c_6,c_7>0$ are constants independent of $n,p,L$.
\end{proposition}

\begin{remark}
Condition \eqref{eq:thml2:moment} is the $L_2$-$L_4$ norm equivalence condition for $\psi(X_t)$; see \cite{mendelson2020robust}. Let $\xi = w^\top \psi(X_t)$. Then it becomes $\EE(\xi^4)\le c_1\big(\EE(\xi^2) \big)^2$, implying that the kurtosis of $\xi$ is bounded. The $L_2$-$L_4$ norm equivalence plays an important role in random matrix theory and it holds in various settings, such as sub-Gaussian random vectors. See \cite{mendelson2020robust} for more details and more examples. 

In addition, to ensure (\ref{eq:eiga}), following \citet{oliveira2016lower}, condition \eqref{eq:thml2:moment} can be relaxed by letting $u$ be sparse vectors satisfying $\| u \|_0\le n$.
\end{remark}

\begin{assumption}\label{asmp:sparsity}
Let $S:=\{(j,k): h_{jk}(\cdot)\not\equiv 0,1\le j, k\le p\}$ and $S_j:=\{k: h_{jk}(\cdot)\not\equiv 0, 1\le k\le p \}$, $1\le j\le p$. Assume that nonzero indices 
\begin{align*}
s_0:=\max_{1\le j\le p}\sum_{k=1}^p \mathbf{1}_{\{h_{jk}\not\equiv 0\}}=\max_{1\le j\le p}\textrm{Card}(S_j) = o(p)\,\,\textrm{and}\,\,
s:=\sum_{j=1}^p\sum_{k=1}^p \mathbf{1}_{\{h_{jk}\not\equiv 0\}} =\textrm{Card}(S)= o(p^2).
\end{align*}
\end{assumption}

Assumption \ref{asmp:sparsity} imposes a sparsity condition on the nonlinear functions. Structural sparsity condition is often used in high-dimensional setting, for example, \cite{cai2011} in covariance matrix estimation. To achieve convergence rates without an additional factor of $p$, as is typically desired in high-dimensional settings, global boundedness of the quantities in Assumption \ref{asmp:sparsity} is usually required, as in \cite{koltchinskii2010sparsity}. However, \cite{raskutti2012minimax} finds an elaborate way to circumvent this requirement when studying sparse additive models with RKHS components.

The following Proposition \ref{prop:bddri} establishes an upper bound on the remainder term $\|r_t\|_\infty$ as a function of the smoothness level $\beta$, the number of basis functions $L$, and the sparsity level $s_0$. Moreover, the quantity $\frac{1}{n}\sum_{t=1}^n [h_{jk}(X_{t-1}^{(k)})- h_{jk}^{(L)}(X_{t-1}^{(k)}) ]^2$ serves as a measure of the $L_2$ bias between $h_{jk}$ and its orthogonal projection onto the finite-dimensional subspace spanned by the chosen basis functions.

\begin{proposition} \label{prop:bddri}
Under Assumptions \ref{asmp:expansion} and \ref{asmp:sparsity}, we have
\begin{align*}
\|r_t\|_\infty =\max_{1\le j\le p} \left| \sum_{k=1}^p [h_{jk}(X_{t-1}^{(k)})- h_{jk}^{(L)}(X_{t-1}^{(k)})] \right| &\le BC(2\beta-1)^{-1} s_0 L^{1/2-\beta} , \\
\max_{1\le j,k\le p} \frac{1}{n}\sum_{t=1}^n \left[h_{jk}(X_{t-1}^{(k)})- h_{jk}^{(L)}(X_{t-1}^{(k)}) \right]^2 &\le B^2 C^2 (2\beta-1)^{-2} L^{1-2\beta} .
\end{align*}
\end{proposition}

Formally, we have the following asymptotic properties for the $\ell_1$ regularized estimators. Theorem \ref{thm:l2} shows how the rate of convergence of $\hat b-b^*$ and the errors of the estimated functions $\hat h_{jk}$ depend on the sparsity of functions, basis expansions, the dependence strength of the processes and the moment condition.

\begin{theorem} \label{thm:l2}
Suppose Assumptions \ref{asmp:LipconstH}, \ref{asmp:moment}(ii), \ref{asmp:expansion},  \ref{asmp:eig} and \ref{asmp:sparsity} hold. 
Let $\hat b$ be the corresponding LASSO solution given in the optimization problem \eqref{problem2}. Consider the estimator 
\begin{eqnarray} \label{h:estm}
\hat h_{jk}(x)= \sum_{l=1}^L \psi_{k,l}(x) \hat b_{j,k}^{(l)}, \qquad 1\le j,k \le p.
\end{eqnarray}
Suppose that condition \eqref{eq:thml2:moment} holds.
Assume that
\begin{align} \label{eq:thml2:lambda}
\lambda \ge c_2 \left(\sqrt{\frac{L \log (pL)}{n}}+s_0L^{1-\beta}  \right),
\end{align}
for some $c_2>0$. Also suppose that 
$$n\ge c_3 s_0L\cdot\log(n)\log(pL)+c_3L^2\cdot\log(n)\log(pL)$$ for some sufficiently large constant $c_3$. We have, with probability approaching one (as $n,p \rightarrow \infty$),
\begin{align} 
\|\hat b-b^*\|_2 &\le c_4\sqrt{s}\lambda, \label{eq:rate:l2a} \\
\sum_{j=1}^p\sum_{k=1}^p \|\hat h_{jk}-h_{jk} \|_2^2 &\le c_5 s\lambda^2 +c_5 sL^{-2\beta}, \label{eq:rate:l2b} \\
\frac1n \sum_{t=1}^n \sum_{j,k=1}^p ( \hat h_{jk} (X_{t-1}^{(k)}) -  h_{jk} (X_{t-1}^{(k)}) )^2 &\le c_6 s \lambda^2 +c_6 sL^{1-2\beta}, \label{eq:rate:l2c}
\end{align}
where $c_4,c_5,c_6>0$ are constants depending on $\rho$ and $\mu_e$.
\end{theorem}

Observe that since $s\le s_0p$, the bounds in \eqref{eq:rate:l2a}, \eqref{eq:rate:l2b} and \eqref{eq:rate:l2c} imply that
\begin{align*} 
\max_{1\le j\le p}\|\hat b_{j}-b_{j}^*\|_2 &\le c_4\sqrt{s_0}\lambda,  \\
\max_{1\le j\le p}\sum_{k=1}^p \|\hat h_{jk}-h_{jk} \|_2^2 &\le c_5 s_0\lambda^2 +c_5 s_0L^{-2\beta}, \\
\frac1n \sum_{t=1}^n \sum_{k=1}^p ( \hat h_{jk} (X_{t-1}^{(k)}) -  h_{jk} (X_{t-1}^{(k)}) )^2 &\le c_6 s_0 \lambda^2 +c_6 s_0 L^{1-2\beta} ,
\end{align*}
where $b^*$ and $b_{j}^*$ are defined in \eqref{eq:def:b} and \eqref{eq:XiPhi}. The quantity $\rho$ measures the dependence strength of the processes, and the constant $\mu_e$ encodes the moment condition. Theorem \ref{thm:l2} shows that, provided Assumptions \ref{asmp:LipconstH} and \ref{asmp:moment}(ii) hold with $\rho\le \rho_0<1$ and $\rho_0$ is a constant, neither the dependence strength nor the moment constant $\mu_e$ affects these convergence rates.
The second terms in \eqref{eq:rate:l2b} and \eqref{eq:rate:l2c} quantify the bias due to truncating the basis expansion. Moreover, Theorem \ref{thm:l2} implies that if the noise $\epsilon_t$ has finite exponential moments, then we may allow the dimension $p$ to grow as fast as $e^{n^c}$ for some constant $0<c<1$; the exponent $c$ depends on the chosen truncation level $L$ of basis expansion.

It is instructive to compare the two terms in the tuning requirement $\lambda$ from \eqref{eq:thml2:lambda}. In the case with relative low dimension $\log(p)\lesssim s_0^2n L^{1-2\beta}$ and low basis number $L\lesssim s_0^{2/(2\beta-1)}(n/\log n)^{1/(2\beta-1)}$, the basis-expansion bias term $s_0L^{1-\beta}$ dominates. On the other hand, if the dimension $p$ is large such that $\log(p)\gtrsim s_0^2n L^{1-2\beta}$ or basis number $L$ is large with $L\gtrsim s_0^{2/(2\beta-1)}(n/\log n)^{1/(2\beta-1)}$, the stochastic term $(n^{-1}L\log(pL))^{1/2}$ becomes the leading factor.

\begin{remark}
The convergence rates of the penalized estimators in \eqref{eq:rate:l2b} and \eqref{eq:rate:l2c} contain two sources of bias: (a) the first from the penalty $\lambda$, and (b) the second from the truncation parameter $L$ (which depends on the smoothness of the function space, $\beta$).
\end{remark}

\begin{remark}[Use of Bernstein-type Inequalities]
Bernstein-type inequalities play a crucial role in the theoretical analysis of high-dimensional methods with regularization. Define the loss function
\begin{align*}
F(b)&=\frac1n \sum_{t=1}^n \| X_t-\Psi(X_{t-1})^\top b \|_2^2 + \lambda\sum_{j,k=1}^p \sqrt{\frac1n \sum_{t=1}^n (\psi_{k} (X_{t-1}^{(k)})^\top  b_{j,k})^2 }, 
\end{align*}
and define
\begin{align*}
\Sigma_{k}&=\frac1n \sum_{t=1}^n \psi_{k} (X_{t-1}^{(k)})\psi_{k} (X_{t-1}^{(k)})^\top   
\quad \mathrm{and}\quad
J_n=\frac1n \sum_{t=1}^n   \Psi(X_{t-1})\Psi(X_{t-1})^\top.  
\end{align*}
Following the standard proof technique for regularized estimators \citep{negahban2012}, we compare $F(\hat b)$ to $F(b^*) $, where $\hat b$ minimizes $F(b)$, to obtain
\begin{align*}
0 \ge F(\hat b)-F(b^*)
=& -2\nabla_n^\top(\hat b-b^*)+(\hat b-b^*)^\top J_n (\hat b-b^*)+\lambda \sum_{j,k=1}^p \big(\|\Sigma_{k}^{1/2}\hat b_{j,k}\|_2
- \|\Sigma_{k}^{1/2} b^*_{j,k}\|_2\big)  ,
\end{align*}
where $\nabla_n$ is the gradient of the least squares loss, defined in \eqref{eq:nablan} below. In our analysis, Theorem \ref{thm:bern} is not applied verbatim in the proof of Theorem \ref{thm:l2}, but its underlying arguments and closely related concentration inequalities are used. 
First, we establish a high probability bound on $|\nabla_n|_{2,\infty}$, where $|\cdot|_{2,\alpha}$ is defined in \eqref{eq:norm} below. In particular, Lemma \ref{lambda:l2:a} requires an exponential-type tail probability bound for $\frac1n\sum_{t=1}^n g(X_{t-1})\epsilon_{t}^{(j)}$ analogous to the bound in Theorem \ref{thm:bern}. Next, we need a high probability bound for the quadratic term $(\hat b-b^*)^\top J_n (\hat b-b^*)$. Obtaining this bound also relies on Bernstein-type inequalities, as generalized in Lemma \ref{lem:mdepinftypart1}.
However, because of temporal dependence, the quantities in Lemmas \ref{lem:mdepinftypart1} and \ref{lambda:l2:a} involve quadratic forms or noise terms rather than simple Lipschitz functions $g$ of $X_t$ as in Assumption \ref{asmp:functiong}, so Theorem \ref{thm:bern} cannot be applied directly. 
We therefore adapt its technical arguments to establish a corresponding exponential-type tail probability bound and then use those bounds to prove Theorem \ref{thm:l2}.

\end{remark}

\begin{remark}
Our framework in Theorem \ref{thm:l2} is quite general: it accommodates a broad class of nonlinear VAR processes whose innovations need not be sub-Gaussian. By contrast, \cite{han2015} and \cite{basu2015} focus on linear VAR models with i.i.d. Gaussian errors, estimating the transition matrix. Like those linear VAR analyses, we also allow the ambient dimension $p$ to vastly exceed the sample size $n$.

A crucial distinction arises in the tuning parameter condition \eqref{eq:thml2:lambda}. The second term on the right, originating from the bias in truncating the basis expansion, enters the gradient of the loss and must be retained when verifying restricted strong convexity \citep{negahban2012}. Consequently, the truncation level $L$ influences both the choice of $\lambda$ and the estimator's convergence rate. 

In the fully nonlinear setting, one typically requires $L\to\infty$, so the first term $\sqrt{L\log (p L)/n}$ in $\lambda$'s bound exceeds the familiar $\sqrt{\log (p)/n}$ rate for linear VARs \citep{basu2015}. This inflation can be viewed as the statistical ``cost of nonlinearity''.
However, in special cases where each $h_{jk}$ admits an exact (or arbitrarily precise) finite dimensional basis representation, the bias term $s_0 L^{1-\beta}$ in \eqref{eq:thml2:lambda} vanishes and the first term collapses to $\sqrt{\log(p)/n}$. Under those circumstances, our nonlinear estimator attains the same tuning and convergence rates as its linear counterpart.

\end{remark}

Next, we turn to model-selection consistency. In place of Assumptions \ref{asmp:eig}, we present an alternative condition that directly targets the support of each component. 
To simplify the notation, let $\Psi_{S_j}(X_{t})=(\psi_{k} (X_{t}^{(k)})^\top, k\in S_j)$ be the truncated feature vector in $\RR^{L\cdot \textrm{Card}(S_j)}$, where $\psi_{k} $ is defined in \eqref{eq:def:b}. We then assemble these vectors into the block-diagonal matrix
\begin{align*}
\Psi_S(X_{t})=\left(\begin{array}{cccccc}
\Psi_{S_1} (X_{t})^\top &  0 & 0 &\cdots  & 0 \\
0 & \Psi_{S_2} (X_{t})^\top  & 0 &\cdots  & 0\\
0 & 0 & \Psi_{S_3} (X_{t})^\top  &\cdots  & 0\\
\vdots & \vdots & \vdots & \ddots & \vdots\\
0 & 0 & 0 &\cdots  & \Psi_{S_p} (X_{t})^\top \\
\end{array}\right) .
\end{align*}

%

\begin{assumption} \label{asmp:eig2}
There are some constants $\phi_{\max}, \phi_{\min}>0, 0< \delta\le 1$, so that 
\begin{align} 
\min_{1\le j\le p}\lambda_{\min}\Big\{\EE \Psi_{S_j}(X_{t-1})\Psi_{S_j}(X_{t-1})^\top\Big\}&\ge \phi_{\min}>0, \label{eq:incoherencea} \\
\max_{1\le j\le p} \lambda_{\max}\Big\{\EE \Psi_{S_j}(X_{t-1})\Psi_{S_j}(X_{t-1})^\top\Big\}&\le \phi_{\max}<\infty, \label{eq:incoherenceb}
\end{align}
and
\begin{align} \label{eq:incoherencec}
\max_{1\leq j\leq p} \left\| \left( \EE \Psi_{S_j^c} (X_{t-1})\Psi_{S_j}(X_{t-1})^\top \right)  \left(\EE \Psi_{S_j}(X_{l-1})\Psi_{S_j}(X_{l-1})^\top \right)^{-1}\right\|_{2,\infty}\le \sqrt{\frac{\phi_{\min}}{\phi_{\max}}} \cdot \frac{1-\delta}{\sqrt{s_0}},
\end{align}
where the induced matrix $(2,\infty)$-norm is defined as $\|A\|_{2,\infty}=\max_{1\le j\le m_1} \sqrt{ \sum_{k=1}^{m_2} A_{jk}^2}$ for $A\in\RR^{m_1\times m_2}$.
\end{assumption}

This assumption corresponds to the condition of \cite{ravikumar2009,ravikumar2010}. Similar to Assumption \ref{asmp:eig}, \eqref{eq:incoherencea} and \eqref{eq:incoherenceb} are also standard, and are commonly imposed for high-dimensional regression analysis. Besides, \eqref{eq:incoherencec} relates to the incoherence condition, see e.g. \cite{wainwright2009,ravikumar2010}. In the following proposition, we establish a sample version of Assumption \ref{asmp:eig2}.

\begin{proposition}\label{prop:incoherence}
Suppose Assumptions \ref{asmp:LipconstH}, \ref{asmp:moment}(ii), \ref{asmp:sparsity} and \ref{asmp:eig2} hold. Assume $\sup_x|\psi_{k,l}(x)|\le B$ for any $1\le k\le p, 1\le l\le L$.  
Assume that \eqref{eq:thml2:sparsity:np} holds and that for some constant $c_1 >0$ does not rely on $p,s_0,L$, such that for all $u\in\mathbb R^{{\rm Card}(S_j) L}$, 
\begin{align}\label{eq:thm3:moment}
\EE(u^\top\Psi_{S_j}(X_t)\Psi_{S_j}(X_t)^\top u)^2\le c_1 \big(u^\top\EE(\Psi_{S_j}(X_t)\Psi_{S_j}(X_t)^\top )u\big)^2.    
\end{align}
Then, with probability approaching one (as $n,p \rightarrow \infty$), we have 
\begin{align} 
\lambda_{\min}\Big\{\frac{1}{n}\sum_{t=1}^n\Psi_S(X_{t-1})\Psi_S(X_{t-1})^\top\Big\}&\ge (1+o(1))\phi_{\min}>0, \label{eq:eig2a} \\
\lambda_{\max}\Big\{\frac{1}{n}\sum_{t=1}^n\Psi_S(X_{t-1})\Psi_S(X_{t-1})^\top\Big\}&\le (1+o(1))\phi_{\max}<\infty, \label{eq:eig2b}
\end{align}
and
\begin{align} \label{eq:eig2c}
\max_{1\leq j\leq p}\max_{k\in S_j^c} \left\| \left( \frac1n  \sum_{t=1}^n\psi_{k} (X_{t-1}^{(k)})\Psi_{S_j}(X_{t-1})^\top \right)  \left(\frac1n  \sum_{l=1}^n\Psi_{S_j}(X_{l-1})\Psi_{S_j}(X_{l-1})^\top \right)^{-1}\right\|_2\nonumber\\
\le (1+o(1))\sqrt{\frac{\phi_{\min}}{\phi_{\max}}} \cdot \frac{1-\delta}{\sqrt{s_0}}.
\end{align}

\end{proposition}

In Theorem \ref{thm:l2:sparsity}, we show that, under certain conditions, our method recovers the sparsity pattern asymptotically. Recall $S=\{(j,k): h_{jk}(\cdot)\not\equiv 0,1\le j, k\le p\}$. Then $S=\{(j,k): b_{j,k}^*\neq 0,1\le j, k\le p\}$. Let $\hat S_n:=\{(j,k): \hat b_{j,k}\neq 0,1\le j, k\le p\}$. 

\begin{theorem} \label{thm:l2:sparsity}
Suppose Assumptions \ref{asmp:LipconstH}, \ref{asmp:moment}(ii), \ref{asmp:expansion}, \ref{asmp:sparsity} and \ref{asmp:eig2} hold. Let $\hat b$ be the corresponding LASSO solution given in the optimization problem \eqref{problem2}. Let $\beta>3/2$. Suppose that condition \eqref{eq:thm3:moment} holds. Assume that
\begin{align} \label{eq:thml2:sparsity:np}
\frac{s_0L^2\cdot\log(pL)}{n}+s_0L^{1-2\beta/3} \rightarrow 0,
\end{align}
and 
\begin{align} \label{eq:thml2:sparsity:lambda}
\lambda \sqrt{s_0} L + \lambda^{-1} \sqrt{\frac{L\log(n)}{n} } + \lambda^{-1}s_0 L^{1-\beta} \rightarrow 0.
\end{align}
Then the solution $\hat b$ to problem \eqref{problem2} is unique and satisfies $\hat S_n=S$, with probability approaching one (as $n,p\to\infty$).
\end{theorem}

In a $p$-dimensional vector time series, the pattern of direct influences among variables can be represented by a binary adjacency matrix $A = (a_{jk}) \in \{0,1\}^{p\times p}$, where
\[
a_{jk} =
\begin{cases}
1, & \text{if variable }k\text{ directly influences variable }j,\\
0, & \text{otherwise.}
\end{cases}
\]
In a linear VAR model $X_t = G\,X_{t-1} + \epsilon_t$, this network structure is typically inferred from the nonzero entries of the transition matrix $G$, which is often assumed to be sparse \citep{hall2018}. A theory-free principle was advocated in \cite{sims1980} for inferring economic relations between variables of linear VARs.

In our nonlinear VAR framework, each component function $h_{jk}$ quantifies the influence of $k$ on $j$. Moreover, the group lasso formulation in \eqref{problem2} yields a sparse estimate $\hat b$, so that many blocks $\hat b_{j,k}$ are exactly zero. We therefore define the estimated adjacency matrix $\hat A=(\hat a_{jk})$ by
\[
\hat a_{jk} =
\begin{cases}
1, & \text{if }\hat b_{j,k}\neq 0,\\
0, & \text{if }\hat b_{j,k}=0.
\end{cases}
\]
Since $\hat A$ need not be symmetric, it encodes a directed graph. Our Theorem \ref{thm:l2:sparsity} then guarantees model selection consistency for $\hat A$, ensuring that the true influence network is recovered with high probability. We demonstrate the proposed network estimation method on real data in Section \ref{section:data}.

\section{Simulation Studies} \label{section:simulation}

In this section, we shall evaluate the numerical performance of the proposed estimation procedures of nonlinear VAR models.

\begin{figure}[htbp!]
\centering
\subfigure[random]{
\begin{minipage}{0.315\linewidth}
\includegraphics[scale=0.46]{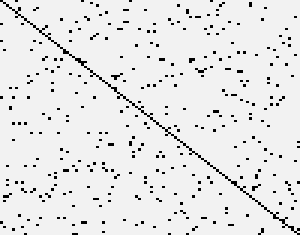}
\end{minipage}}
\subfigure[band]{
\begin{minipage}{0.315\linewidth}
\includegraphics[scale=0.46]{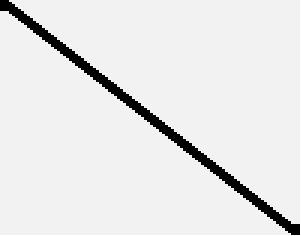}
\end{minipage}}
\subfigure[cluster]{
\begin{minipage}{0.315\linewidth}
\includegraphics[scale=0.46]{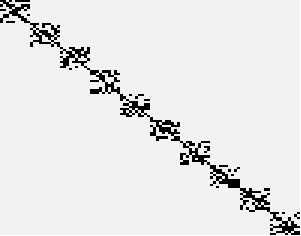}
\end{minipage}}
\caption{Three different network matrix patterns used in the simulation studies. Here gray points represent the zero entries and black points represent nonzero entries.}
\label{figure:pattern}
\end{figure}

We design three different patterns of the binary transition matrix (network matrix, see Section \ref{section:model}) $A$: random, band, cluster. Typical realizations of these patterns are illustrated in Figure \ref{figure:pattern}. The pattern ``cluster'' has block diagonal structure, where each block is of dimension $10\times 10$ and satisfies the pattern ``random''. In each dimension $j$, $1\le j\le p$, we randomly assign 5 nonzero functions, according to the pattern of the transition matrix. The relevant nonzero component functions are given by
\begin{align*}
f_1(x)&=0.2x, \\
f_2(x)&=-0.15\sin(1.5x), \\
f_3(x)&=-0.5\Phi(x,0.5,1), \\
f_4(x)&=0.2xe^{-0.5x^2}, \\
f_5(x)&=0.15\log(|x|+2),
\end{align*}
where $\Phi(\cdot,0.5,1)$ is the Gaussian probability distribution function with mean 0.5 and standard deviation 1. In other words, for each $j$ with $1\le j\le p$, we randomly select 5 functions $h_{jk} $ ($1\le k\le p$) to be the above nonzero functions. The rest $p-5$ functions of $h_{jk}$ ($1\le k\le p$) are all zeros. Elementary calculation shows that this nonlinear VAR process is stable and satisfies Assumption \ref{asmp:LipconstH}. In order to ensure reasonable signal to noise ratio, the error processes $\epsilon_t$ are generated from $0.2 N(0,1)$.

In all the conducted experiments, we assess the model selection performance of our procedure using the area under the receiver operating characteristic curve (AUROC) and the area under the Precison-Recall curve (AUPR) ignoring the sign (positive negative influence), where the ROC curve is created by plotting the true positive rate (TPR) against the false positive rate (FPR) and the precision-recall curve is a plot of the precision against the recall. Define TPR, FPR, precision and recall as follows
\begin{align*}
\rm
TPR=recall=\frac{TP}{TP+FN}, \quad FPR=\frac{FP}{TN+FP}, \quad Precision= \frac{TP}{TP+FP}.  
\end{align*}
Here TP and TN stand for true positives and true negatives, respectively, and FP and FN stand for false positives/negatives. We choose a set of data dimensions $p=20, 50,100$ while the sample size is $n= 50,100,200,500$, respectively. The empirical values reported in Tables \ref{tab:table1} are averages over 1000 replications.

It can be seen from Table \ref{tab:table1} that the proposed estimation procedure of nonlinear VAR model performs fairly well as reflected in both AUROC and AUPR. In particular, when the sample size is moderate ($n\ge 100$), our method provides pretty good AUROC in all cases. As expected, when the sample size $n$ increases, our method performs better. And both AUROC and AUPR decreases as the dimension $p$ increase. Besides, our proposed method makes no significant differences in terms of 3 patterns of transition matrix.

\begin{table}[htbp]
\centering
\caption{Model selection performance of the proposed nonlinear VAR method with three different patterns of the transition matrix, ``random'', ``band'', ``cluster'', based on 1000 replications.}
\label{tab:table1}
\begin{tabular}{cccccccccccccccccccccccccc} \hline 
$p$& & \multicolumn{4}{c}{AUROC} & & \multicolumn{4}{c}{AUPR}  \\ \cline{3-6}\cline{8-11}
& $n$ & 50 & 100 & 200 & 500 & & 50 & 100 & 200 & 500\\ \hline
& \multicolumn{5}{l}{Pattern ``random''}\\ \hline
20   && 0.633 & 0.744 & 0.851 & 0.924 & & 0.443 & 0.651 & 0.856 & 0.937 \\
50   && 0.611 & 0.720 & 0.842 & 0.920 & & 0.230 & 0.458 & 0.753 & 0.904 \\
100   && 0.591 & 0.696 & 0.830 & 0.918 & & 0.132 & 0.320 & 0.666 & 0.883 \\ \hline
& \multicolumn{5}{l}{Pattern ``band''}\\ \hline
20   && 0.647 & 0.753 & 0.858 & 0.928 & & 0.469 & 0.681 & 0.864 & 0.938 \\
50   && 0.610 & 0.720 & 0.841 & 0.920 & & 0.234 & 0.464 & 0.758 & 0.905 \\
100   && 0.592 & 0.698 & 0.830 & 0.918 & & 0.143 & 0.339 & 0.672 & 0.881 \\\hline
& \multicolumn{5}{l}{Pattern ``cluster''}\\ \hline
20   && 0.642 & 0.746 & 0.855 & 0.922 & & 0.464 & 0.667 & 0.861 & 0.933 \\
50   && 0.609 & 0.718 & 0.839 & 0.920 & & 0.231 & 0.454 & 0.744 & 0.905 \\
100   && 0.591 & 0.696 & 0.827 & 0.918 & & 0.138 & 0.328 & 0.661 & 0.883 \\
\hline
\end{tabular}
\end{table}

\section{Real Data Analysis} \label{section:data}
We now apply our nonlinear VAR model to the analysis of a real biological gene regulatory network time series expression data. The network is an {\it E. coli} SOS DNA repair system, which has been well studied in biology, see e.g, \cite{ronen2002}. The main function of the SOS signaling pathway is to regulate cellular immunity and repair DNA damage. We consider an eight gene network, part of the SOS DNA repair network in the bacteria {\it E. coli}. The time series gene expression data set of the network was collected by \cite{ronen2002}. The data are kinetics of 8 genes, that is, lexA, recA, ruvA, polB, umuDC, uvrA, uvrD, uvrY, where lexA and recA are the key genes in the pathway. The 8 genes were measured at 50 instants which are evenly spaced by 6 min intervals.

\begin{figure}[htbp!]
\centering
\includegraphics[scale=0.38]{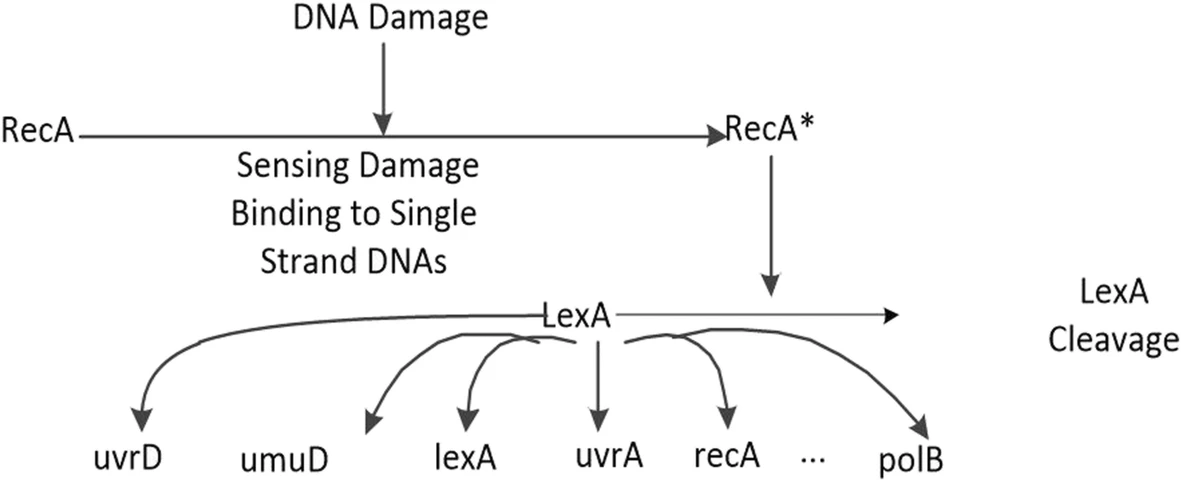}
\caption{The bacterial SOS DNA repair system}
\label{figure:sos}
\end{figure}

\begin{figure}[htbp!]
\centering
\includegraphics[scale=0.9]{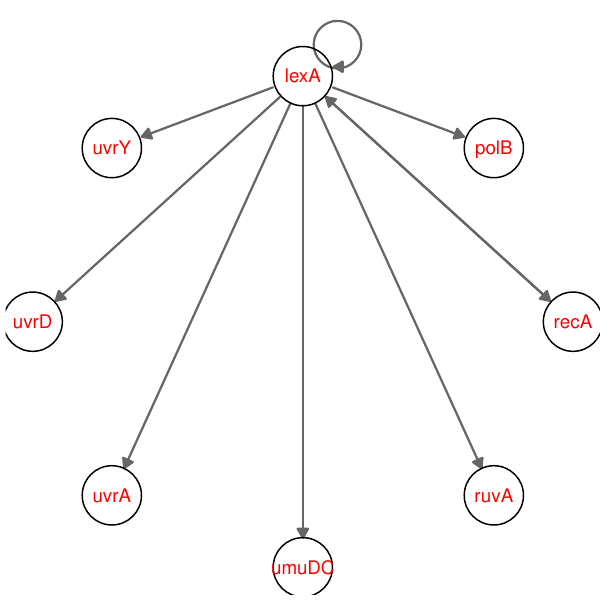}
\caption{The target SOS DNA repair network}
\label{figure:sos_true}
\end{figure}

\begin{figure}[htbp!]
\centering
\includegraphics[scale=0.9]{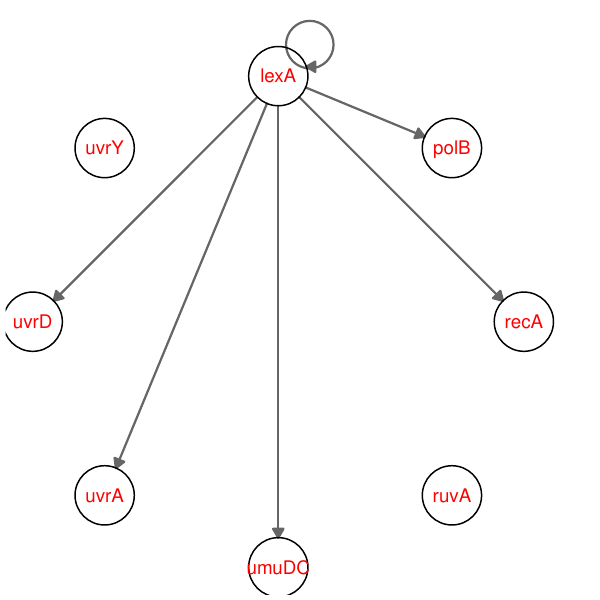}
\caption{Reconstruction of SOS DNA repair network by nonlinear VAR model}
\label{figure:sos_non}
\end{figure}

\begin{figure}[htbp!]
\centering
\includegraphics[scale=0.9]{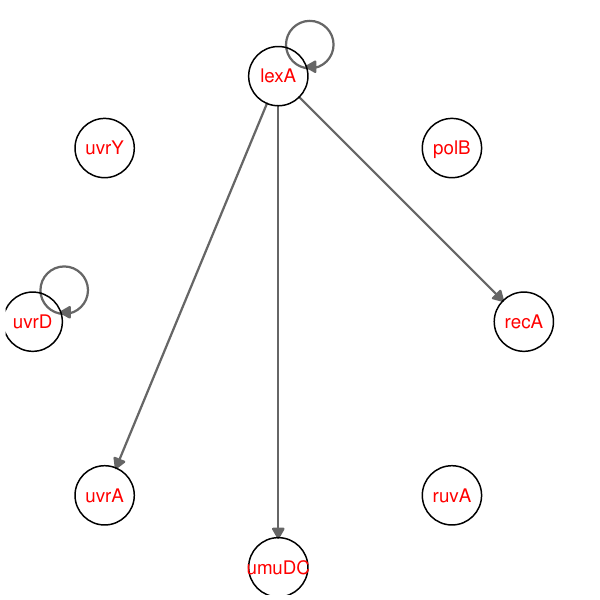}
\caption{Reconstruction of SOS DNA repair network by linear VAR model}
\label{figure:sos_linear}
\end{figure}

We compare the performance of our method with the Lasso regularized linear VAR method (\cite{basu2015}). The tuning parameter $\lambda$ in both methods and the number of basis function $L$ are chosen by time series cross-validation procedure (see \cite{han2015}). Figure \ref{figure:sos} represents the bacterial SOS DNA repair system. Figure \ref{figure:sos_true} shows the real SOS DNA repair network, which contains 9 edges. Figures \ref{figure:sos_non} and \ref{figure:sos_linear} show the inferred gene regulatory networks using our nonlinear VAR model and the $\ell_1$ regularized linear VAR model, respectively. In Figure \ref{figure:sos_non}, one can see that our method finds 6 out of the 9 edges in the target network and identifies lexA as the hub gene for this network. Our method identifies most interactions except lexA$\to$ruvA, lexA$\to$uvrY and recA$\to$lexA. In comparison, in the Figure \ref{figure:sos_linear}, the $\ell_1$ regularized linear VAR model recognizes only 4 out of the 9 true edges, and predicts a wrong edge. Furthermore, our proposed method gives the area under ROC curve 0.8116 and the area under Precison-Recall curve 0.6836. While, the $\ell_1$ regularized linear VAR model gives AUROC 0.7222 and AUPR 0.6036. In summary, our proposed method has a better performance than the regularized linear VAR model on the SOS DNA repair network, although none of these two methods can faithfully recover all of the edges. This phenomenon also confirms that there exists nonlinear dynamics in the gene regulatory networks.

\section{Proofs} \label{section:proofs}
Write random variable $\bm\xi \in {\mathcal L}^m $, $m\ge 1$, if the $m$-norm $\|\bm \xi \|_m := (\EE |\bm\xi|^m)^{1/m} < \infty$. Denote $\|\bm \xi \|:=\|\bm \xi \|_2$. 
Let $\F_k=(\ldots,\epsilon_{k-1},\epsilon_k)$, $\F_k^n=\{\epsilon_{k},\ldots,\epsilon_n\}$, and $\EE_0(X)=X-\EE X.$ Define projection operator $P_k(\cdot)=\EE(\cdot|\F_k)-\EE(\cdot|\F_{k-1})$, $k\in\ZZ$. Let $(\epsilon_k')_{ k\in\ZZ}$ be an i.i.d. copy of $(\epsilon_k)_{k\in\ZZ}$, so that $\epsilon_i,$ $\epsilon_j'$, $i,j\in\ZZ$ are i.i.d. For any $X_t=\mathcal G(\ldots,\epsilon_{t-1},\epsilon_t)$, where $\mathcal G$ is a measurable function, we define the coupled version $X_{t,\{k\}}=\mathcal G(\ldots,\epsilon_{k-1},\epsilon_{k}',\epsilon_{k+1},\ldots,\epsilon_t)$. If $k>i$, then $X_{t,\{k\}}=X_t$.

\subsection{Proofs of Theorems in Section \ref{section:bern}}

\begin{lemma}[\cite{burkholder1988}, \cite{rio2009}]
\label{lem:burkholder}
Let $q>1$, $q'=\min\{q,2\}$. Let $D_T=\sum_{t=1}^T \xi_t,$ where $\xi_t\in \L^q$ are martingale differences. Then
\begin{align*}
\|D_T\|_q^{q'} \leq K_q^{q'}\sum_{t=1}^T\|\xi_t\|_q^{q'}, \mbox{ where } K_q = \max\big\{(q-1)^{-1}, \sqrt{q-1}\big\}.
\end{align*}
\end{lemma}

\begin{lemma}
\label{lem:increaseorder}
Let $\epsilon\in\RR^p$ be a random vector with non-negative entries, satisfying Assumption \ref{asmp:moment}(i) with $\mu_q<\infty$ for some $q\geq 2$. For non-negative vectors $v_t\in\mathbb R^p$, $i\geq 0,$ assume $\|v_t\|_1\leq \rho^t$ where $\rho<1.$ Denote
\begin{align*}
X=\sum_{t=0}^\infty \min\big\{  v_t^\top \epsilon,M\big\}.   
\end{align*}
Take $c_0=-\rho^2\log\rho/(2e)$. Then for any $c\leq c_0/M$, $\EE(e^{cX})$ exists and
\begin{align*}
\EE (e^{c_0 X/M})-\EE (c_0X/M)-1\leq  \mu_2^2M^{-2}<\infty.
\end{align*}

\end{lemma}

\begin{proof}
Note that we have the decomposition
\begin{align*}
X   
=M \sum_{t=0}^\infty  \1_{\{v_t^\top \epsilon\ge  M\}}
+\sum_{t=0}^\infty  v_t^\top \epsilon\1_{ \{v_t^\top \epsilon<  M\}}=:\I_1+\I_2.
\end{align*}
For $\I_1$ part, we have for any $m\geq 1,$
\begin{align}
\label{eq:bddi1mdraft}
\EE|\I_1|^m 
&\le  M^m\left(\sum_{t=0}^\infty\big\|\1_{\{v_t^\top \epsilon\ge  M\}}\big\|_m\right)^m
=
M^m\Big(\sum_{t=0}^\infty \PP( v_t^\top \epsilon\geq  M)^{1/m}\Big)^m.
\end{align}
By Markov's inequality,
\begin{align*}
  \PP(  v_t^\top \epsilon\geq  M)\le  \|v_t^\top\epsilon\|_2^2/M^2
  \le \rho^{2i}\mu_2^2/M^2.
\end{align*}
Applying above into \eqref{eq:bddi1mdraft}, we further have
\begin{align*}
\EE|\I_1|^m \le M^m \Big( \mu_2^{2/m} M^{-2/m} \sum_{t=0}^\infty \rho^{2i/m} \Big)^m 
\le \mu_2^2(1-\rho^{2/m})^{-m} M^{m-2}.    
\end{align*}
Since for any $m\geq 1$,
we have
\begin{align} \label{eq:rhobdlg}
1-\rho^{2/m}\geq (1-\rho^2)/m.
\end{align}
We further obtain 
\begin{align*}
\EE|\I_1|^m &
\leq \mu_2^2((1-\rho^2)/m)^{-m} M^{m-2}.
\end{align*}
Choose $c_{1,M}=-\rho^2\log(\rho)/(eM),$ then by $m!\geq (2\pi)^{1/2}m^{m+1/2}e^{-m}$ (\cite{robbins1955}),  we have
\begin{align*}
\sum_{m\geq 2}\frac{\EE((c_{1,M}\I_1)^m)}{m!}\leq \frac12 \mu_2^2 M^{-2}.
\end{align*}
For $\I_2$ part, for any $m\geq 2,$ 
\begin{align*}
\EE |\I_2|^m
&\leq \left(\sum_{t=0}^\infty \| v_t^\top \epsilon\1_{ v_t^\top \epsilon <  M }\|_m\right)^m
\leq \left(\sum_{t=0}^\infty (M^{m-2}\EE| v_t^\top \epsilon|^2)^{1/m}\right)^m
\leq \mu_2^2 \Big(M^{1-2/m}\sum_{t=0}^\infty \rho^{tq/m}\Big)^m \\
&\leq  \mu_2^2(-2\rho^2 \log(\rho)/m)^{-m}M^{m-2} ,
\end{align*}
where the last inequality is by \eqref{eq:rhobdlg}. Therefore
\begin{align*}
\sum_{m\geq 2}\frac{\EE((c_{1,M}\I_2)^m)}{m!}\leq \frac12 \mu_2^2M^{-2}<\infty,
\end{align*}
We complete the proof by combining the two parts and setting $c_0=Mc_{1,M}/2$.
\begin{align*}
\EE e^{c_0X/M}-1-\EE(c_0X/M)=\sum_{m\geq 2} \frac{\EE((c_0X/M)^m)}{m!}
\leq \sum_{m\geq 2} \frac{\EE((c_{1,M}\I_1)^m)}{m!}+\sum_{m\geq 2} \frac{\EE((c_{1,M}\I_2)^m)}{m!}
\leq \mu_2^2M^{-2}.
\end{align*}

\end{proof}


\begin{proof}[\bf Proof of Theorem \ref{thm:bern}]{\bf Part (i).}
Without loss of generality, assume $\|G\|_1=1$. For $X_t=\mathcal G(\ldots,\epsilon_{t-1},\epsilon_t)$, where $\mathcal G$ is some measurable function, we define the coupled version $$X_{t,\{k\}}=\mathcal G(\ldots,\epsilon_{k-1},\epsilon_{k}',\epsilon_{k+1},\ldots,\epsilon_t).$$ Recall that $\abs(x)=(|x_1|,\ldots,|x_p|)^\top$ for $x=(x_1,\ldots,x_p)^\top$. By Assumption \ref{asmp:LipconstH}, for $k\leq i-1,$ we have
\begin{align*}
\abs(X_t-X_{t,\{k\}})
\leq H \abs(X_{t-1}-X_{t-1,\{k\}}),
\end{align*}
and for $k=i$, $\abs(X_t-X_{t,\{k\}})=\epsilon_k-\epsilon_k'.$
Hence by induction, we obtain 
\begin{align*}
\abs(X_t-X_{t,\{k\}})
\leq H^{t-k} \abs(\epsilon_k-\epsilon_k').
\end{align*}
Since the function 
$g$ is Lipschitz continuous, combined with the above inequality, we have
\begin{align}
\label{eq:pkgxi}
|P_k g(X_t)|
&=\big|\EE(g(X_t)-g(X_{t,\{k\}})|\F_k)\big|\nonumber\\
&\leq \EE \left( G^\top \abs(X_t-X_{t,\{k\}})  \Big|\F_k\right)\nonumber\\
&\leq  \EE \left( G^\top  H^{t-k} \abs(\epsilon_k-\epsilon_k') \Big|\F_k\right).
\end{align}
Let $S_n(g)=\sum_{t=1}^n \big(g(X_t)-\EE g(X_t)\big).$ For $k\leq n,$ denote $\xi_k=P_k(S_n(g)).$ Then $$S_n(g)=\sum_{k\leq n}\xi_k.$$
The tail probability can be decomposed into two parts
\begin{align*}
\PP(S_n(g)\geq 2z)\leq \PP\Big(\sum_{-n<k\leq n}\xi_k\geq z\Big)+\PP\Big(\sum_{k\leq -n}\xi_k\geq z\Big)=:\I_{1}+\I_{2}.    
\end{align*}
In the following, we will first bound $\xi_k$ and then address $\I_1$ and $\I_2$ separately. The first part $\I_1$ is the leading term, while the second part $\I_2$ is relatively small.
By Assumption \ref{asmp:LipconstH} and $\|G\|_1\leq 1$, we have
\begin{align*}
\|H^{t-k \top}G\|_1\leq \|H\|_\infty^{t-k}\|G\|_1\leq \rho^{t-k}.    
\end{align*}
 Denote $v_t=H^{t\top}G$. Since $|g|_\infty\leq M$, we have $|P_k g(X_t)|\leq 2M.$ Thus by \eqref{eq:pkgxi}, 
\begin{align}
\label{eq:usexik1}
|\xi_k|\leq \sum_{t=1}^n |P_kg(X_t)|\leq \sum_{t=k\vee 1}^n\min\Big\{v_{t-k}^\top\EE(\abs(\epsilon_k-\epsilon_k') \big |\F_k), 2M \Big\},\quad\textrm{with } \|v_t\|_1\leq \rho^t.  
\end{align}

For part $\I_{1},$ let $h^*:=-\rho^2 (\log\rho)/(4eM).$ By Lemma \ref{lem:increaseorder} and \eqref{eq:usexik1} for any $0<h\leq h^*$, 
$\EE(e^{|\xi_k| h})<\infty.$
Note that $\EE(\xi_k|\F_{k-1})=0.$ Then 
\begin{align}
\label{eq:uligt_1}
\EE(e^{\xi_kh}|\F_{k-1})
&=1+\EE(e^{\xi_kh}-\xi_kh-1|\F_{k-1})\nonumber\\
&\leq 1+\EE\Big[\frac{e^{|\xi_k|h}-|\xi_k|h-1}{h^2}\Big|\F_{k-1}\Big]h^2,
\end{align}
in view of $e^x-x\leq e^{|x|}-|x|$ for any $x.$ Note that for any fixed $x>0,$ $(e^{tx}-tx-1)/t^2$ is increasing in $t \in (0,\infty)$. By Lemma \ref{lem:increaseorder}, we have
\begin{align}
\label{eq:uligt_2}
\EE\Big[\frac{e^{|\xi_k|h}-|\xi_k|h-1}{h^2}\Big|\F_{k-1}\Big]
&\leq \EE\Big[\frac{e^{|\xi_k|h^*}-|\xi_k|h^*-1}{h^{*2}}\Big|\F_{k-1}\Big]\nonumber\\
&\leq (h^*)^{-2}\mu_2^2 (2M)^{-2}\nonumber\\
&\leq c_3,
\end{align}
where $c_3=4e^2(-\rho^2 \log \rho)^{-2}\mu_2^2$. Hence for any $h\leq h^*,$ by \eqref{eq:uligt_1} and \eqref{eq:uligt_2},
\begin{align}
\label{eq:cqgammaref}
\EE(e^{\xi_kh}|\F_{k-1})\leq 1+c_3 h^2.
\end{align}
By Markov's inequality we have 
$\I_1\leq e^{-zh}\EE[\exp(\sum_{-n <k\leq n} \xi_k h) ].$ Let $h=\min\{z(4c_3n)^{-1}, h^*\},$ then by recursively applying \eqref{eq:cqgammaref},
\begin{align}
\I_1
\le & e^{-zh}\EE\Big(e^{ \sum_{k=-n+1}^{n-1} \xi_k h}\EE(e^{\xi_nh}|\F_{n-1}) \Big)\nonumber\\
\leq & e^{-zh}(1+c_3h^2)^{2n}\nonumber\\
\le & \exp\big(-zh+2n c_3 h^2\big) \nonumber\\
\leq & \exp\Big\{ - \frac{z^2}{8c_3n+ c_4Mz}\Big\},\label{eq:fistnpart}
\end{align} 
where the third inequality is due to $1+x \leq e^{x}$ for $x>0,$ and $c_4=8e/(-\rho^2 \log\rho).$

For $\I_2,$ by \eqref{eq:usexik1}, $\|\xi_k\|_q\leq \sum_{t=1}^n\rho^{t-k}\mu_q\leq \rho^{1-k}(1-\rho)^{-1}\mu_q,$ for $k\leq 0.$ Then by Lemma \ref{lem:burkholder},
\begin{align} \label{eq:secondnpart}
\I_2 & \leq z^{-q}\Big((q-1)\sum_{k\leq -n} \|\xi_k\|_q^2\Big)^{q/2} \notag\\
&\le (q-1)^{q/2}z^{-q}\Big(\sum_{k\leq -n} \|\xi_k\|_q^2\Big)^{q/2}  \notag \\
&\leq c_5 \rho^{qn}/z^q=c_5 e^{-qn\log(\rho^{-1})}/z^q,
\end{align}
where $c_5=(q-1)^{q/2}\mu_q^q(1-\rho)^{-3q/2}$ only depends on $\rho, q$ and $\mu_q$.

Combining $\I_1$ and $\I_2$ parts, the desired result follows by noticing $z\leq 2Mn$.
\ \\


\noindent {\bf Part (ii).}
Without loss of generality, assume $\|G\|_1=1$. Similar to the proof of Theorem \ref{thm:bern}(i), let $S_n(g)=\sum_{t=1}^n \big(g(X_t)-\EE g(X_t)\big)$, and $\xi_k=P_k(S_n(g)).$ Then $S_n(g)=\sum_{k\leq n}\xi_k$, and 
\begin{align*}
\PP(S_n(g)\geq 2z)\leq \PP\Big(\sum_{-n<k\leq n}\xi_k\geq z\Big)+\PP\Big(\sum_{k\leq -n}\xi_k\geq z\Big)=:\I_{1}+\I_{2}.    
\end{align*}
Denote $v_t=H^{t\top}G$ and $\omega_k=\sum_{t=1 \vee k}^n v_{t-k}$. Since \eqref{eq:pkgxi} still holds, we have 
\begin{align} \label{eq:usexik2}
|\xi_k|\leq \sum_{t=k\vee 1}^n v_{t-k}^\top\EE(\abs(\epsilon_k-\epsilon_k') \big |\F_k)=\omega_k^\top \EE(\abs(\epsilon_k-\epsilon_k') \big |\F_k).
\end{align}
For $\I_{2},$ $k\leq -n,$ $\|w_k\|_1\leq \rho^{1-k}/(1-\rho).$ Let $h^*:=c_0(1-\rho)/\rho.$ By \eqref{eq:uligt_1} and \eqref{eq:uligt_2}, for any $0\le h\le h^*$,
\begin{align}
\label{eq:pdfse1}
\EE(e^{\xi_kh}|\F_{k-1})
&\leq 1+\EE\Big[\frac{e^{|\xi_k|h^*}-|\xi_k|h^*-1}{h^{*2}}\Big|\F_{k-1}\Big]h^2
\leq  1+\frac{\EE( e^{|\xi_k|h^*}-1|\F_{k-1})}{h^{*2}} h^2.
\end{align}
Let $a_k=\rho^{1-k}/(1-\rho)$ and $u_k=w_k/a_k,$ then
\begin{align*}
\EE( e^{|\xi_k|h^*}-1|\F_{k-1}) 
&\leq \EE\Big( e^{w_k^\top \abs(\epsilon_k-\epsilon_k') h^* } -1\Big) = \EE\Big( e^{c_0u_k^\top \abs(\epsilon_k-\epsilon_k') \rho^{-k}} -1 \Big).
\end{align*}
If $f(0)=0$, then $\EE(f(X))=\int_0^\infty f'(t)\PP(X\geq t)\d t.$ Therefore we further obtain
\begin{align}
\label{eq:pdfse2}
\EE( e^{|\xi_k|h^*}-1|\F_{k-1})   
&\leq \int_0^\infty e^{t\rho^{-k}}\rho^{-k}\PP(c_0u_k^\top \abs(\epsilon_k-\epsilon_k') \geq t)\d t\nonumber\\
&\leq \rho^{-k} \int_0^\infty e^{-t(1-\rho^{-k})}\mu_e^2\d t
\leq \rho^{-k} (1-\rho)^{-1}\mu_e^2.
\end{align}
Since $1+x\leq e^x,$ by \eqref{eq:pdfse1} and \eqref{eq:pdfse2}, 
\begin{align}
\label{eq:cqgammaref2}
\EE(e^{\xi_kh}|\F_{k-1})\leq 1+ \rho^{-k} (1-\rho)^{-1}\mu_e^2(h^{*})^{-2} h^2\leq e^{c_3\rho^{-k}h^2},    
\end{align}
where $c_3=\mu_e^2(1-\rho)^{-3}\rho^2c_0^{-2}.$
Recursively applying \eqref{eq:cqgammaref2}, we can obtain
\begin{align*}
\I_2\leq
 e^{-zh^*}\EE\Big( e^{\sum_{k\leq -n}\xi_kh^*}\Big)
 \leq \exp(-zh^*+c_4\rho^n h^{*2} ),
\end{align*}
where $c_4=c_3/(1-\rho).$
Similar to \eqref{eq:fistnpart}, we can bound the $\I_1$ part and we complete the proof.
\end{proof}

\subsection{Proofs of Theorems in Section \ref{section:var}}

By \eqref{eq:def:b}, for vector $b=(b_{j,k})_{1\leq j,k\leq p}$ and $b_{j,k}\in\RR^L$, define the $(2,\alpha)$ group structure norm 
\begin{align} \label{eq:norm}
 |b|_{2,\alpha}:=||b_{j,k}|_2  |_{\alpha} = \left(\sum_{j=1}^p\sum_{k=1}^p \left(\sum_{l=1}^L (b_{j,k}^{(l)})^2\right)^{\alpha/2} \right)^{1/\alpha}, 
\end{align}
where $\alpha\ge 1$. For instance, with the choice $\alpha=1$, this norm corresponds to the regularizer that underlies the group Lasso. For $\alpha=\infty$,
\begin{align*}
 |b|_{2,\infty}:=||b_{j,k}|_2  |_{\infty} = \max_{1\le j,k\le p} \left(\sum_{l=1}^L (b_{j,k}^{(l)})^2\right)^{1/2}. 
\end{align*}

\begin{proof}[\bf Proof of Proposition \ref{prop:bddri}]
Note that since basis functions are orthonormal, $\|h_{jk}\|_2=(\sum_{l=1}^\infty (b_{jl}^{(l)*})^2)^{1/2}.$ 
Since basis functions are bounded by $B$, by Assumption \ref{asmp:expansion}, we have
\begin{align*}
\|h_{jk}-h_{jk}^{(L)}\|_\infty
&\le \sum_{l\ge L+1}|b_{jk}^{(l)*}|B \\
&= B \sum_{l\ge L+1}\frac{ |b_{jk}^{(l)*}| l^\beta} {l^\beta} \\ 
&\le B \sqrt{ \sum_{l\ge L+1} (b_{jl}^{(l)*})^2 l^{2\beta}} \sqrt{\sum_{l\ge L+1} l^{-2\beta}} \\
&\le BC(2\beta-1)^{-1}L^{1/2-\beta} .
\end{align*}
Hence, as $s_0=\max_{1\le j\le p}$Card$(S_j)$ with $S_j:=\{k: h_{jk}(\cdot) \neq 0, 1\le k\le p \}$,
\begin{align*}
\|r_t\|_\infty \le \sum_{k=1}^p \|h_{jk}-h_{jk}^{(L)}\|_\infty
\le BC(2\beta-1)^{-1} s_0 L^{1/2-\beta} .
\end{align*}
Furthermore, we have
\begin{align*}
\frac{1}{n}\sum_{t=1}^n \left[h_{jk}(X_{t-1}^{(k)})- h_{jk}^{(L)}(X_{t-1}^{(k)}) \right]^2 &= \frac{1}{n}\sum_{t=1}^n \left[\sum_{l\ge L+1} \psi_{k,l}(X_{t-1}^{(k)}) b_{jk}^{(l)*}  \right]^2 \\
&\le B^2 \left[\sum_{l\ge L+1} b_{jk}^{(l)*}  \right]^2 \\
&\le B^2 C^2 (2\beta-1)^{-2} L^{1-2\beta} .
\end{align*}
Then we obtain the desired result. 
\end{proof}

\begin{proof}[\bf Proof of Proposition \ref{prop:eig}]
We first prove part (i). By \eqref{eq:eig1}, we have, for any $u\in\mathbb R^{pL}$ with $\|u\|_2=1$,
\begin{eqnarray*}
 \EE  u^\top \psi (X_t) \psi (X_t)^\top u  \geq \phi_{L}.
\end{eqnarray*}
Let $m=4(-\log \rho)^{-1}\log(n)$. Recall $\F_k^n=\{\epsilon_{k},\ldots,\epsilon_n\}.$ By Lemma \ref{lemma:rscddiffpartmd1}, we have, with probability at least $1-mp^{-c_1}/12- 2mp L e^{-3n/(10m)}$, for any $u\in\mathbb R^{pL}$, 
\begin{align*}
&\frac{1}{n}\sum_{t=1}^nu^\top\EE\big( \psi (X_t)\psi (X_t)^\top |\F_{t-m+1}^n\big) u \ge \frac12 u^\top\EE \psi (X_t)\psi (X_t)^\top u-\frac{c_2\log(n)\log(pL)}{n} \left\|u \right\|_1^2.
\end{align*}
Note that $L=o(n)$. Let $z=1$ in Lemma \ref{lem:mdepinftypart}, we can obtain, with probability at least $1-mp^{-c_1}/12- 2mp L e^{-3n/(10m)}-e^{-c_3 n}$, for any $u\in\mathbb R^{pL}$, 
\begin{align*}
& \frac{1}{n}\sum_{t=1}^nu^\top\big( \psi (X_t)\psi (X_t)^\top \big) u \ge \frac12 u^\top\EE \psi (X_t)\psi (X_t)^\top u-\frac{c_2\log(n)\log(pL)}{n} \left\|u \right\|_1^2-\frac1n \|u\|_2^2.
\end{align*}
Then \eqref{eq:eiga} follows.

For part (ii), denote $\Omega_{k}=\EE(\psi_{k} (X_{t}^{(k)}) \psi_{k} (X_{t}^{(k)})^\top)$. For $m=o(n)$, let $N=[(n-1)/m]$ and $\mathcal N=\{1,m+1,2m+1,\ldots,(N-1)m+1\}$. Then there exists constant $c_3>0$ such that for any $1\leq l_1,l_2\leq L,$ $z>0$, we have
\begin{align*}
\PP\left( \Big|\frac{1}{N}\sum_{t\in\N}\EE\bigg( (\psi_{k} (X_{t}^{(k)}) \psi_{k} (X_{t}^{(k)})^\top)_{l_1,l_2}|\F_{t-m+1}^n\Big)-
(\Omega_{k})_{l_1,l_2}\Big|
\geq z \right)
\leq 2\exp\big\{-c_3Nz^2\big\}.
\end{align*}
Therefore with probability at least $1-2L^2\exp\{-c_3Nz^2\},$ for any $u\in\mathbb R^{L}$ with $\|u\|_2=1$,
\begin{align*}
\Big|\frac{1}{N}\sum_{t\in\N}\EE\bigg( u^\top \psi_{k} (X_{t}^{(k)}) \psi_{k} (X_{t}^{(k)})^\top u|\F_{t-m+1}^n\Big)-
u^\top \Omega_{k}u\Big|
\leq Lz.
\end{align*}
Take $z=c_4\sqrt{ (\log(p)+\log(L))/N} $ for some constant $c_4$ large enough. Then we have with probability greater than 
$1-m(pL)^{-c_4}$, for any $u\in\mathbb R^L$, $\|u\|_2=1$, $1\le k\le p$,
\begin{eqnarray*}
\frac{1}{n}\sum_{t=1}^n\EE\bigg( u^\top \psi_{k} (X_{t}^{(k)}) \psi_{k} (X_{t}^{(k)})^\top u|\F_{t-m+1}^n\Big)
\leq  \phi_{U} + c_{5} L\sqrt{\frac{\log(p)+\log(L)}{N}}.
\end{eqnarray*}
Then \eqref{eq:eigb} follows by combining above and Lemma \ref{lem:mdepinftypart} with $z=1$ and 
$m=4(-\log \rho)^{-1}\log(n)$.
\end{proof}

\begin{lemma}\label{lemma:rscddiffpartmd1}
For $m=o(n)$, denote $N=[(n-1)/m]$ and $\mathcal N=\{1,m+1,2m+1,\ldots,(N-1)m+1\}$.
Consider the VAR process \eqref{eq:autoregmodel}, suppose Assumptions \ref{asmp:LipconstH} and \ref{asmp:moment}(ii) hold. Assume that there exists a constant $c>0$, such that for all $u\in\mathbb R^{p L}$, $\EE[(u^\top\psi(X_t)\psi(X_t)^\top u)^2]\le c\big(u^\top\EE(\psi(X_t)\psi(X_t)^\top )u\big)^2$. Let $N\ge C \log(p L)$, where $C>0$ is a sufficiently large constant. Then, we have, with probability at least $1-p^{-c_1}/12- 2p L e^{-3N/10}$, 
\begin{align*}
&\forall u\in \RR^{p L}, \frac{1}{N}\sum_{t\in\mathcal N}u^\top\EE\big( \psi(X_t)\psi(X_t)^\top |\F_{t-m+1}^n\big) u \ge \frac12 u^\top\EE \psi(X_t)\psi(X_t)^\top u-\frac{c_2\log(p L)}{N} \left\|u \right\|_1^2, 
\end{align*}
where $c_1>0$ is a sufficiently large constant and $c_2$ depends only on $c$ and $B$.
\end{lemma}

\begin{proof}
Recall for any $1\le k\le p, 1\le l\le L$, $\sup_x|\psi_{k,l}(x)|\le B$, some $B\geq 1$, and $\F_k^n=\{\epsilon_{k},\ldots,\epsilon_n\}.$ 
Denote $\Sigma=\EE (\psi(X_t)\psi(X_t)^\top) $ and $$\tilde\Sigma_N=N^{-1}\sum_{t\in\mathcal N}\EE\big( \psi(X_t)\psi(X_t)^\top |\F_{t-m+1}^n\big).$$ Let $\tilde\Sigma_{\text{diag}}$ be the diagonal of $\tilde\Sigma_N$. Note that $\EE\big( \psi(X_t)\psi(X_t)^\top |\F_{t-m+1}^n\big)=\EE\big( \psi(X_t)\psi(X_t)^\top |\F_{t-m+1}^t\big)$ are independent for all $t\in\mathcal N$. By Jensen's inequality, 
$$\EE\left[\big(\EE \big(u^\top\psi(X_t)\psi(X_t)^\top u |\F_{t-m+1}^n \big)\big)^2 \right] \le \EE[(u^\top\psi(X_t)\psi(X_t)^\top u)^2] \le c\big(u^\top\EE(\psi(X_t)\psi(X_t)^\top )u\big)^2.$$
Then, employing similar arguments as in the proof of Lemmas 5.1 and 5.2 in \cite{oliveira2016lower}, we can obtain, for $N\ge 1568 c(c_3+1)\log(p L)$ and $c_3>0$, 
\begin{eqnarray} \label{eq:eigen1a}
\PP\left(\forall u\in \RR^{p L}, u^\top\tilde\Sigma_N u \ge \frac12 u^\top\Sigma u-\frac{1568c(c_3+1)\log(p L)}{N} \left| \tilde\Sigma_{\text{diag}}^{1/2}u \right|_1^2  \right)\ge 1-\frac{1}{12}p^{-c_3} .
\end{eqnarray}
Since for any $1\le k\le p, 1\le l\le L$, $|\psi_{k,l}|_\infty\le B$, then, by Bernstein's inequality, we have,
\begin{eqnarray*}
\PP\left( \left| \frac{1}{N}\sum_{t\in\mathcal N} (\psi_{k,l}(X_{t}^{(k)})^2-\EE [\psi_{k,l}(X_{t}^{(k)})^2 |\F_{t-m+1}^n ]) \right|\ge z \right) \le 2\exp\left(-\frac{Nz^2}{2B^4+4B^2z/3}\right).
\end{eqnarray*}
Hence, we have
\begin{eqnarray*}
\PP\left( \max_{1\le k\le p,1\le l\le L} \left| \frac{1}{N}\sum_{t\in\mathcal N} \psi_{k,l}(X_{t}^{(k)})^2 \right|\ge 2B^2 \right) \le 2p L\exp\left(-10N/3 \right).
\end{eqnarray*}
Combining the above inequality with \eqref{eq:eigen1a}, it follows that, with probability at least $1-p^{-c_3}/12- 2p L e^{-3N/10}$, for any $u\in \RR^{p L}$,
\begin{eqnarray*}
u^\top\tilde\Sigma_N u \ge \frac12 u^\top\Sigma u-\frac{3136B^2c(c_3+1)\log(p L)}{N} \left\|u \right\|_1^2  . 
\end{eqnarray*}
\end{proof}

\begin{lemma}(m-approximation)
\label{lem:mdepinftypart}
Considering the VAR process \eqref{eq:autoregmodel}, suppose Assumptions \ref{asmp:LipconstH} and \ref{asmp:moment} (ii) hold. Let $z\rho^{-m}/(s_0L) >C n$, where $C>0$ is a sufficient large constant. We have
\begin{align*}
&\PP\left(\sup_{\|u\|_2=1,\ \|u\|_1^2= s_0L}\Big|\sum_{t=1}^n u^\top \big[\psi (X_t) \psi (X_t)^\top-\EE\big(\psi (X_t) \psi (X_t)^\top |\F_{t-m+1}^n\big)\big]u\Big|\ge z\right)
\leq   s_0^2L^2 e^{-cn},
\end{align*}
for some constant $c>0$.
\end{lemma}
\begin{proof}
For matrix $A,$ denote by $A_{k_1,k_2}$ the $(k_1,k_2)th$ entry of $A$, and let
$\EE_{t-m+1}(\cdot)=(\cdot)-\EE(\cdot|\F_{t-m+1}^n),$ then we have
\begin{align*}
&\PP\left(\sup_{\|u\|_2=1,\ \|u\|_1^2= s_0L}\Big|u^\top\sum_{t=1}^n\EE_{t-m+1}\big(\psi (X_t) \psi (X_t)^\top \big)u\Big|\ge z\right)    \\
\leq &\PP\left(\sup_{\|u\|_2=1,\ \|u\|_1^2=s_0L}\|u\|_1^2\max_{1\leq k_1,k_2\leq p L}\Big|\sum_{t=1}^n\EE_{t-m+1}\big((\psi (X_t) \psi (X_t)^\top)_{k_1,k_2}\big)\1_{u_{k_1},u_{k_2}\neq 0}\Big|\ge z\right)  \\
\leq &  s_0^2L^2\max_{1\leq k_1,k_2\leq p L}\PP\left(\Big|\sum_{t=1}^n\EE_{t-m+1}\big((\psi (X_t) \psi (X_t)^\top)_{k_1,k_2}\big)\Big|\ge z/(s_0L)\right).  
\end{align*}
By construction, for any indices $k_1,k_2,$ there exist functions 
$$\phi_1,\phi_2\in \{f : \RR^p \to \RR | f(x)=\psi_{k,l}(x_k) \textrm{ for some } 1\leq k\leq p,1\leq l\leq L \}$$
such that $(\psi (X_t) \psi (X_t)^\top)_{k_1,k_2}=\phi_1(X_t)\phi_2(X_t).$ 
Since function $\psi_{k,l}$ satisfies conditions in Lemma \ref{lem:mdepinftypart1}, we complete the proof. 
\end{proof}

\begin{lemma}\label{lem:mdepinftypart1}
Consider the VAR process \eqref{eq:autoregmodel}, suppose Assumptions \ref{asmp:LipconstH} and \ref{asmp:moment}(ii) hold. Assume functions $\phi_1, \phi_2:\RR^p\rightarrow \RR$ are both bounded with $|\phi_t|_\infty\leq B,$ $i=1,2.$ For any $x,y\in\RR^p,$ assume
$|\phi_t(x)-\phi_t(y)|\leq \beta^\top \abs(x-y)=\sum_{j=1}^p \beta_j|x_j-y_j|,$ where $\|\beta\|_1\leq 1.$ Then we have
\begin{align}
\PP\Big(\big|\sum_{t=1}^n \big[\phi_1(X_t)\phi_2(X_t)-\EE\big(\phi_1(X_{t})\phi_2(X_{t})|\F_{t-m+1}^n\big)\big]  \big|\geq z\Big)  
\leq e^{-c\min\{n, z\rho^{-m}, z^2\rho^{-2m}/n\}},
\end{align}
where constant $c$ only depends on $\rho,$ $\mu_2$, $\mu_e$ and $B$.  
\end{lemma}

\begin{proof}
Recall $\F_k^n=\{\epsilon_{k},\ldots,\epsilon_n\}.$ Denote
\begin{align*}
S_n=\sum_{t=1}^n \big[\phi_1(X_t)\phi_2(X_t)-\EE\big(\phi_1(X_{t})\phi_2(X_{t})|\F_{t-m+1}^n\big)\big]
\quad\textrm{and}\quad \xi_k=\EE(S_n|\F_{k-1}^n)-\EE(S_n|\F_{k}^n).
\end{align*} 
Then 
$S_n=\sum_{k\leq n-m+1}\xi_k$ and
\begin{align}
\label{eq:xik}
|\xi_k|
&\leq \sum_{t=(k+m-1)\vee 1}^n \EE\Big( |\phi_1(X_{t,\{k\}})-\phi_1(X_t)||\phi_2(X_t)| \big|\F_k^n\Big)\nonumber\\
&\quad +\sum_{t=(k+m-1)\vee 1}^n \EE\Big( |\phi_1(X_{t,\{k\}})||\phi_2(X_{t,\{k\}})-\phi_2(X_{t})|\big|\F_k^n\Big)
=:\xi_{1k}+\xi_{2k}.
\end{align}
Since $|\phi_1(X_{t,\{k\}})-\phi_1(X_t)|\leq \beta^\top H^{t-k}\abs(\epsilon_k'-\epsilon_k)$ and $|\phi_1|_\infty\leq B$, we have
\begin{align*}
\xi_{1k}\le \sum_{t=(k+m-1)\vee 1}^n B\cdot \EE\Big( \beta^\top H^{t-k} \abs(\epsilon_k'-\epsilon_k) \big|\F_k^n\Big).
\end{align*}
A similar bound can be derived for $\xi_{2k}.$ Hence
\begin{align*}
|\xi_k|\le \EE\big(\omega_k^\top \abs(\epsilon_k'-\epsilon_k) |\F_k^n\big),   \mbox{ where } 
 \omega_k^\top=2B\beta^\top \sum_{t=(k+m-1)\vee 1}^n  H^{t-k}.
\end{align*}
Then $\|\omega_k\|_1\leq 2B(1-\rho)^{-1}\rho^{m-1}$ for $k> -n$ and 
$\|\omega_k\|_1\leq 2B(1-\rho)^{-1}\rho^{1-k}$ if $k\leq -n.$
For $k\leq -n,$ since $\xi_k$ are martingale differences, by Burkholder's inequality (Lemma \ref{lem:burkholder}), we have, for any $q\ge 2$,
\begin{align*}
\Big\|\sum_{k\leq -n}\xi_k\Big\|_q^2
\le (q-1)^{q/2} \Big(\sum_{k\le -n} \|\xi_k\|_q^2 \Big)^{q/2} 
\le (q-1)^{q/2} (2B)^q \mu_q^q (1-\rho)^{-q} (1-\rho^2)^{-q/2} \rho^q \rho^{nq}.
\end{align*}
Thus by Markov's inequality
\begin{align*}
\PP\Big(\big|\sum_{k\le -n}\xi_k\big|\geq z\Big)\le
z^{-2}4B^2(1-\rho)^{-2}(1-\rho^2)^{-1}\mu_2^2\rho^2\cdot\rho^{2n}
\le z^{-2}4B^2(1-\rho)^{-4}\mu_2^2\rho^2\cdot e^{-(-2\log \rho)n}.
\end{align*}
For $k>-n,$ let $h^*=(2B)^{-1}(1-\rho)\rho c_0$ and $\xi_k'=\xi_k/\rho^m$ Then $\EE \exp(h^*|\xi_k'|)\leq 2\mu_e<\infty.$
By \eqref{eq:uligt_1}, \eqref{eq:uligt_2} and \eqref{eq:cqgammaref}, we have for any $h\leq h^*$,
\begin{align*}
\EE(e^{\xi_k'h}|\F_{k-1})\leq 1+c_1 h^2,  
\end{align*}
where $c_1=2\mu_e h^{*-2}.$
Similar as \eqref{eq:fistnpart}, we have
\begin{align*}
\PP\Big(\big|\sum_{k=-n+1}^n\xi_k/\rho^m\big|\geq z\Big)
\leq \inf_{h\leq h^*}\exp\big(-z h+2c_1 n h^{2}\big) \leq \exp\big\{- z^2/(c_2z+c_3n)\big\},
\end{align*} 
for some constants $c_2, c_3$ depending on $\rho,\mu_2,\mu_e$ and $B$. Then the desired result follows.
\end{proof}

\begin{remark}
The proof of Lemma \ref{lem:mdepinftypart1} follows a similar approach to that of Theorem \ref{thm:bern}.  
\end{remark}

\begin{proof}[\bf Proof of Theorem \ref{thm:l2}]
Let $$F(b)=\frac1n \sum_{t=1}^n \| X_t-\Psi(X_{t-1})^\top b \|_2^2 + \lambda\sum_{j,k=1}^p \sqrt{\frac1n \sum_{t=1}^n (\psi_{k} (X_{t-1}^{(k)})^\top  b_{j,k})^2 }.$$
Define 
\begin{align}
\label{eq:nablan}
\nabla_n= \frac1n \sum_{t=1}^n \Psi(X_{t-1})( X_t-\Psi(X_{t-1})^\top b^* ).    
\end{align}
Recall the definition of $|\cdot|_{2,\alpha}$ in \eqref{eq:norm}. Then
\begin{align}
\label{eq:bddnablan}
|\nabla_n|_{2,\infty}
&= \big|\frac1n \sum_{t=1}^n \Psi(X_{t-1})(\epsilon_t+r_t)\big|_{2,\infty} \nonumber \\
&\leq \frac1n \sum_{t=1}^n L^{1/2}\big\|\Psi(X_{t-1})\big\|_{\infty}\|r_t\|_{\infty}+ \big|\frac1n \sum_{t=1}^n \Psi(X_{t-1})\epsilon_t\big|_{2,\infty} \nonumber\\
&:=\I_1+\I_2.
\end{align}
For $\I_1$ part, by \eqref{eq:XiPhi} and Proposition \ref{prop:bddri}, we have $\|\Psi(X_{t-1})\|_{\infty}\leq B$ and thus
$\I_1\leq B^2C (2\beta-1)^{-1} s_0L^{1-\beta}.$ 
For $\I_2$ part, by Lemma \ref{lambda:l2:a}, with probability at least $1-(pL)^{-c'}$, 
$\I_2\leq c\sqrt{L\log(pL)/n},$ for some constants $c,c'>0.$

For $c_2\geq 12(c+ CB^2(2\beta-1)^{-1})/\phi_L$, by Proposition \ref{prop:eig}, we have 
\begin{align*}
\lambda\geq (12/\phi_L)\big(c\sqrt{L\log(pL)/n}+ B^2C(2\beta-1)^{-1} s_0L^{1-\beta} \big)\geq 12|\nabla_n|_{2,\infty}/\phi_L.    
\end{align*}
Let
\begin{align*}
\widetilde\phi_L = \frac{\phi_L}{2}-\frac1n-\frac{c_4(s_0L)\log(n)\log(pL)}{n},
\end{align*} 
and
\begin{align*}
\widetilde\phi_U = \phi_U+c_7L\sqrt{\frac{\log(n)\log(pL)}{n}},
\end{align*}
where $\|u\|_1=s_0L$ in Proposition \ref{prop:eig}, and $c_4,c_7$ are the constants in \eqref{eq:eiga} and \eqref{eq:eigb}. 
Then, for $n\ge c_3 (s_0 L) \log(n)\log(pL)+c_3L^2\log(n)\log(pL)$ with sufficient large constant $c_3>0$, we have 
\begin{align*}
\widetilde\phi_L \ge \frac{\phi_L}{3} \text{ and } \widetilde\phi_U \le 2\phi_U . 
\end{align*}
Denote
\begin{align*}
\Sigma_{k}=\frac1n \sum_{t=1}^n \psi_{k} (X_{t-1}^{(k)})\psi_{k} (X_{t-1}^{(k)})^\top   
\quad \mathrm{and}\quad
J_n=\frac1n \sum_{t=1}^n   \Psi(X_{t-1})\Psi(X_{t-1})^\top.  
\end{align*}
Hence, by Assumption \ref{asmp:eig} and Proposition \ref{prop:eig}, with probability approaching one, we have
\begin{align*} 
F(b)-F(b^*)
=& -2\nabla_n^\top(b-b^*)+(b-b^*)^\top J_n (b-b^*)+\lambda \sum_{j,k=1}^p \big(\|\Sigma_{k}^{1/2}b_{j,k}\|_2
- \|\Sigma_{k}^{1/2} b^*_{j,k}\|_2\big)   \\
\geq &-2|\nabla_n|_{2,\infty}|b-b^*|_{2,1}+\widetilde \phi_L\|b-b^*\|_2^2+\lambda\sum_{j,k\notin S}\|\Sigma_{k}^{1/2}b_{j,k}\|_2-\lambda\sum_{j,k\in S}\|\Sigma_{k}^{1/2}(b_{j,k}-b^*_{j,k})\|_2     \\
\geq & \widetilde\phi_L \|b-b^*\|_2^2 - \lambda(\phi_L/6+\widetilde\phi_U) \sum_{j,k\in S}\|b_{j,k}-b^*_{j,k}\|_2 \\
\geq & (\phi_L/3) \|b-b^*\|_2^2 - \lambda(\phi_L/6+2\phi_U) \sum_{j,k\in S}\|b_{j,k}-b^*_{j,k}\|_2.
\end{align*}
Since Card$(S)=|S|_0=s,$ we have
\begin{align*}
\sum_{j,k\in S}\|b_{j,k}-b^*_{j,k}\|_2 \le \sqrt{s} \sqrt{\sum_{j,k\in S}\|b_{j,k}-b_{j,k}^*\|_2^2} \le s^{1/2}\|b-b^*\|_2.
\end{align*}
Hence $\|\hat b-b^*\|_2\leq (1/2+6\phi_U/\phi_L)\sqrt{s}\lambda$ in view of $F(\hat b)-F(b^*)\leq 0.$

Furthermore,
\begin{align*}
\sum_{j,k=1}^p \|\hat h_{jk}-h_{jk}\|_2^2 &\le \sqrt{2} \sum_{j,k=1}^p \left\| \sum_{l=1}^L(\hat b_{j,k}^{(l)}- b_{j,k}^{(l)*})\psi_{k,l} \right\|_2^2 + \sqrt{2} \sum_{j,k=1}^p \left\| \sum_{l=L+1}^\infty  b_{j,k}^{(l)*} \psi_{k,l} \right\|_2^2. 
\end{align*}
Since $(\psi_{k,l})_{j,k,l}$ are orthonormal basis functions, we have
\begin{align*}
\sum_{j,k=1}^p \|\hat h_{jk}-h_{jk}\|_2^2
&\le \sqrt{2} \sum_{j,k=1}^p \sum_{l=1}^L(\hat b_{j,k}^{(l)}- b_{j,k}^{(l)*})^2 + \sqrt{2} \sum_{j,k=1}^p \sum_{l=L+1}^\infty (b_{j,k}^{(l)*})^2 \\
&\lesssim  s\lambda^2 + \sum_{j,k=1}^p \sum_{l=L+1}^\infty (b_{j,k}^{(l)*})^2 l^{2\beta} l^{-2\beta} \\
&\lesssim  s\lambda^2 +s L^{-2\beta},
\end{align*}
which also implies \eqref{eq:rate:l2b}.

Moreover,
\begin{align*}
\frac1n \sum_{t=1}^n \sum_{j,k=1}^p ( \hat h_{jk} (X_{t-1}^{(k)}) -  h_{jk}^{(L)} (X_{t-1}^{(k)}) )^2 &= (\hat b-b^*)^\top J_n (\hat b-b^*)\\
&\lesssim \sum_{j,k=1}^p (\hat b-b^*)^2 = \| \hat b - b^*\|_2^2 \\
&\lesssim s \lambda^2.
\end{align*}
By Proposition \ref{prop:bddri}, we can obtain \eqref{eq:rate:l2c}.
\end{proof}

\begin{lemma}\label{lambda:l2:a}
For function $g:\RR^p\rightarrow \RR$, assume $|g|_\infty\leq B$. Under Assumption \ref{asmp:moment}(ii), we have
\begin{align}
\PP\Big(\big| \frac1n\sum_{t=1}^n g(X_{t-1})\epsilon_{t}^{(j)}\big|\geq z\Big)\leq  
\begin{cases}
2\exp\Big(-\frac{nz^2}{4c_1}\Big), &\textrm{ if } z\leq 2c_0c_1B^{-1},\\
2\exp\big(-c_0nz/(2B)\big),  &\textrm{ if }  z>2c_0c_1B^{-1}, 
\end{cases}  
\end{align}
where $c_1=\mu_e c_0^{-2} B^2$.
\end{lemma}

\begin{proof}
Let $\xi_t=g(X_{t-1})\epsilon_{t}^{(j)}$. Then $\xi_t, 1\leq i\leq n,$ are martingale differences with respect to $\F_t.$ Let $h^*=c_0/B$. By Assumption \ref{asmp:moment} (ii), for any $0<h\leq h^*$, $\EE(e^{|\xi_k| h})<\infty.$
Since $\EE(\xi_k|\F_{k-1})=0$ and $e^x-x\leq e^{|x|}-|x|$ for any $x,$ we have
\begin{align}
\label{eq:uligt_a}
\EE(e^{\xi_kh}|\F_{k-1})
&=1+\EE(e^{\xi_kh}-\xi_kh-1|\F_{k-1})\nonumber\\
&\leq 1+\EE\Big[\frac{e^{|\xi_k|h}-|\xi_k|h-1}{h^2}\Big|\F_{k-1}\Big]h^2.
\end{align}
Note that for any fixed $x>0,$ $(e^{tx}-tx-1)/t^2$ is increasing in $t \in (0,\infty)$. Hence
\begin{align}
\label{eq:uligt_b}
\EE\Big[\frac{e^{|\xi_k|h}-|\xi_k|h-1}{h^2}\Big|\F_{k-1}\Big]
\leq \EE\Big[\frac{e^{|\xi_k|h^*}-|\xi_k|h^*-1}{h^{*2}}\Big|\F_{k-1}\Big]
\leq \frac{\EE (e^{Bh^*|\epsilon_{t}^{(j)}|})}{h^{*2}}
\leq c_1,
\end{align}
where $c_1=\mu_eB^2c_0^{-2}.$ Combining \eqref{eq:uligt_a} and \eqref{eq:uligt_b}, we can obtain
\begin{align*}
\EE(e^{\xi_kh}|\F_{k-1})\leq 1+c_1h^2.
\end{align*}
Then, by recursively applying the above inequality, we have
\begin{align*}
\PP\Big(\frac1n\sum_{t=1}^n \xi_t \geq z\Big)
&\le  e^{-nzh}\EE\Big(e^{ \sum_{t=1}^{n-1} \xi_t h}\EE(e^{\xi_n h}|\F_{n-1}) \Big) \\
&\le e^{-nzh} (1+c_1 h^2)^{n}  \\
&\le \exp\big(-nzh+n c_1 h^2\big). 
\end{align*}
Take $h=\min\{h^*, z/(2c_1)\},$ we further obtain
\begin{align*}
\PP\Big(\frac1n\sum_{t=1}^n \xi_t \geq z\Big)
\leq & \exp\Big(-\frac{nz^2}{4c_1}\Big)\1_{\{h^*\geq z/(2c_1)\}}+ \exp\big(-c_0nz/(2B)\big)\1_{\{h^*< z/(2c_1)\}}.
\end{align*} 
Similar argument can be applied to $\PP(n^{-1}\sum_{t=1}^n \xi_t \leq -z)$ and the desired result follows.
\end{proof}

\begin{remark}
The proof of Lemma \ref{lambda:l2:a} follows a similar approach to that of Theorem \ref{thm:bern}.  
\end{remark}

\begin{proof}[\bf Proof of Proposition \ref{prop:incoherence}]
Note that 
\begin{align*} 
\lambda_{\min}\Big\{\frac{1}{n}\sum_{t=1}^n\Psi_S(X_{t-1})\Psi_S(X_{t-1})^\top\Big\}&= \min_{1\le j\le p}\lambda_{\min}\Big\{\frac{1}{n}\sum_{t=1}^n\Psi_{S_j}(X_{t-1})\Psi_{S_j}(X_{t-1})^\top\Big\}, \\
\lambda_{\max}\Big\{\frac{1}{n}\sum_{t=1}^n\Psi_S(X_{t-1})\Psi_S(X_{t-1})^\top\Big\}&= \max_{1\le j\le p} \lambda_{\max}\Big\{\frac{1}{n}\sum_{t=1}^n\Psi_{S_j}(X_{t-1})\Psi_{S_j}(X_{t-1})^\top\Big\}.
\end{align*}    
Then, under \eqref{eq:thml2:sparsity:np}, applying \eqref{eq:thm3:moment} and similar arguments in the proofs of Proposition \ref{prop:eig}, we have, in an event $\Omega_1$ with probability approaching one (as $n,p \rightarrow \infty$), 
\begin{align*} 
\min_{1\le j\le p}\lambda_{\min}\Big\{\frac{1}{n}\sum_{t=1}^n\Psi_{S_j}(X_{t-1})\Psi_{S_j}(X_{t-1})^\top\Big\}& \ge (1+o(1))\phi_{\min}>0, \\
\max_{1\le j\le p} \lambda_{\max}\Big\{\frac{1}{n}\sum_{t=1}^n\Psi_{S_j}(X_{t-1})\Psi_{S_j}(X_{t-1})^\top\Big\} & \le (1+o(1))\phi_{\max}>0 .
\end{align*} 
Thus, in the event $\Omega_1$, \eqref{eq:eig2a} and \eqref{eq:eig2b} hold.

Define
\begin{align*}
\widehat Q_{S_j,S_j}&=  \frac{1}{n}\sum_{t=1}^n\Psi_{S_j}(X_{t-1})\Psi_{S_j}(X_{t-1})^\top, \\
\widehat Q_{S_j^c,S_j}&=  \frac{1}{n}\sum_{t=1}^n\Psi_{S_j^c}(X_{t-1})\Psi_{S_j}(X_{t-1})^\top, \\
Q_{S_j,S_j}&=  \EE \Psi_{S_j}(X_{t-1})\Psi_{S_j}(X_{t-1})^\top, \\
Q_{S_j^c,S_j}&=  \EE \Psi_{S_j^c}(X_{t-1})\Psi_{S_j}(X_{t-1})^\top.
\end{align*}
Then, similar to \cite{ravikumar2010}, we decompose the sample matrix as follows
\begin{align*}
\widehat Q_{S_j^c,S_j} \widehat Q_{S_j,S_j}^{-1} &=  Q_{S_j^c,S_j} (\widehat Q_{S_j,S_j}^{-1} - Q_{S_j,S_j}^{-1})
+(\widehat Q_{S_j^c,S_j}-Q_{S_j^c,S_j}) Q_{S_j,S_j}^{-1} \\
&\quad +(\widehat Q_{S_j^c,S_j}-Q_{S_j^c,S_j}) (\widehat Q_{S_j,S_j}^{-1} - Q_{S_j,S_j}^{-1})
+ Q_{S_j^c,S_j}  Q_{S_j,S_j}^{-1} \\
&= \I_1 + \I_2 + \I_3 + \I_4.
\end{align*}
Similar to the proofs of Proposition \ref{prop:eig}, Lemma \ref{lem:mdepinftypart} and Lemma \ref{lem:mdepinftypart1}, we can show in an event $\Omega_{2}$ with probability approaching one (as $n,p \rightarrow \infty$),
\begin{align*}
\| \widehat Q_{S_j^c,S_j}-Q_{S_j^c,S_j} \|_{2,\infty}  =o(1).  
\end{align*}
Based on the properties of the induced matrix norms ($\|\cdot\|_{2,2}=\|\cdot \|_2$), we have in the event $\Omega_1$,
\begin{align*}
\| I_1\|_{2,\infty}   &\le  \| Q_{S_j^c,S_j} Q_{S_j,S_j}^{-1}\|_{2,\infty} \|Q_{S_j,S_j} (\widehat Q_{S_j,S_j}^{-1} - Q_{S_j,S_j}^{-1}) \|_{2,2} \\
& \le \| Q_{S_j^c,S_j} Q_{S_j,S_j}^{-1}\|_{2,\infty} \|Q_{S_j,S_j} \|_2 \|\widehat Q_{S_j,S_j}^{-1} - Q_{S_j,S_j}^{-1} \|_{2} \\
& \le o(1) \| Q_{S_j^c,S_j} Q_{S_j,S_j}^{-1}\|_{2,\infty} .
\end{align*}
Similarly, in the event $\Omega_1\cap\Omega_{2}$
\begin{align*}
\| I_2\|_{2,\infty}   &\le  \| \widehat Q_{S_j^c,S_j}-Q_{S_j^c,S_j} \|_{2,\infty} \|Q_{S_j,S_j}^{-1} \|_{2,2} =o(1), \\
\| I_3\|_{2,\infty}   &\le  \| \widehat Q_{S_j^c,S_j}-Q_{S_j^c,S_j} \|_{2,\infty} \|\widehat Q_{S_j,S_j}^{-1} - Q_{S_j,S_j}^{-1} \|_{2,2} =o(1) .
\end{align*}
It follows that in the event $\Omega_1\cap\Omega_{2}$,
\begin{align*}
\| \widehat Q_{S_j^c,S_j} \widehat Q_{S_j,S_j}^{-1} \|_{2,\infty} &\le (1+o(1) ) \| Q_{S_j^c,S_j} Q_{S_j,S_j}^{-1}\|_{2,\infty}  + o(1).  
\end{align*}
Thus, in the event $\Omega_1\cap\Omega_{2}$ with probability approaching one (as $n,p \rightarrow \infty$), \eqref{eq:eig2a}, \eqref{eq:eig2b} and \eqref{eq:eig2c} hold.
\end{proof}

\begin{proof}[\bf Proof of Theorem \ref{thm:l2:sparsity}]
Let $b_S=(b_{j,k}, (j,k)\in S)\in\RR^{sL}$, and
\begin{align*}
\Omega(b)=\sum_{j,k=1}^p \sqrt{\frac1n \sum_{t=1}^n (\psi_{k} (X_{t-1}^{(k)})^\top b_{j,k})^2 }.
\end{align*}
Denote
\begin{align*}
\hat\Sigma_{S,S}= \frac1n \sum_{t=1}^n \Psi_S(X_{t-1}) \Psi_S(X_{t-1})^\top,
\end{align*}
and
\begin{align*}
\hat\Sigma_{S_j,S_j}= \frac1n \sum_{t=1}^n \Psi_{S_j}(X_{t-1}) \Psi_{S_j} (X_{t-1})^\top.
\end{align*}
By Assumption \ref{asmp:eig2} and Proposition \ref{prop:incoherence}, \eqref{eq:eig2a}, \eqref{eq:eig2b} and \eqref{eq:eig2c} hold on some event $\mathcal Z$ with $\PP(\mathcal Z)\rightarrow 1$. In the following, we shall only work on $\mathcal Z$.

A vector $\hat b\in\RR^{p^2L}$ is an optimum of the objective function in \eqref{problem} if and only if there is a subgradient $\hat g\in \partial \Omega(\hat b)$, such that 
\begin{equation} \label{eq:eq}
\frac2n \sum_{t=1}^n \Psi (X_{t-1}) (\Psi (X_{t-1})^\top \hat b - X_t) +\lambda \hat g =0.
\end{equation}
The subdifferential $\partial \Omega(b)$ is the set of vectors $g=(g_{jk}, 1\le j,k\le p)$, with $\hat g_{jk}\in \RR^{L}$, satisfying
\begin{align}
& g_{jk} = \frac{\frac1n \sum_{t=1}^n \psi_{k} (X_{t-1}^{(k)}) \psi_{k} (X_{t-1}^{(k)})^\top b_{j,k}  }{\sqrt{\frac1n \sum_{t=1}^n (\psi_{k} (X_{t-1}^{(k)})^\top b_{j,k})^2 }},  \label{eq:kkt1}\\
& g_{jk}^\top \left( \frac1n\sum_{t=1}^n \psi_{k} (X_{t-1}^{(k)}) \psi_{k} (X_{t-1}^{(k)})^\top \right)^{-1} g_{jk} \le 1. \label{eq:kkt2}
\end{align}
Following the primal dual witness argument in \cite{ravikumar2009} and \cite{wainwright2009}, it suffices to set $\hat b_{S^c}=0$ and $\hat g_S=\partial \Omega(b^*)_S$, and then show
\begin{align}
&\hat b_{j,k} \neq 0 , \quad \text{ for }(j,k) \in S, \label{eq:signal} \\
&\hat g_{jk}^\top \left( \frac1n\sum_{t=1}^n \psi_{k} (X_{t-1}^{(k)}) \psi_{k} (X_{t-1}^{(k)})^\top \right)^{-1} \hat g_{jk} < 1, \quad \text{ for }(j,k) \in S^c, \label{eq:subgradient}
\end{align}
hold with probability approaching 1.

\begin{enumerate}
\item[(i).] Proof of \eqref{eq:signal}.
\end{enumerate}

Since $\hat b_{S^c}=b_{S^c}^*=0$, \eqref{eq:eq} reduces to
\begin{equation} 
\label{eq:usedlatersoon}
\frac2n \sum_{t=1}^n \Psi_S (X_{t-1}) (\Psi_S (X_{t-1})^\top \hat b_S - X_t) +\lambda \hat g_S =0.
\end{equation}
It implies that
\begin{align}\label{eq:b:decom}
\hat b_S-b_S^* =\hat \Sigma_{S,S}^{-1} \cdot \frac1n\sum_{t=1}^n\Psi_S(X_{t-1}) \epsilon_t + \hat \Sigma_{S,S}^{-1} \cdot \frac1n\sum_{t=1}^n\Psi_S(X_{t-1}) r_t - \frac{\lambda}{2}\hat \Sigma_{S,S}^{-1} \cdot \hat g_S := \I_1+\I_2-\I_3.
\end{align}
We now proceed to bound $\I_1,\I_2$ and $\I_3$. Recall the definition of $|\cdot|_{2,\alpha}$ in \eqref{eq:norm}. Also recall that $\|A\|_{\infty}$ is the matrix $\infty$ norm of $A=(a_{ij})_{n\times m}$ with $\|A\|_{\infty}=\max_{1\le i\le n}\sum_{j=1}^m|a_{ij}|$.

For $\I_1$, we have
\begin{align*}
|\I_1|_{2,\infty} &\le \sqrt{L} \left\| \hat \Sigma_{S,S}^{-1} \right\|_\infty \cdot \left\| \frac1n\sum_{t=1}^n\Psi_S(X_{t-1}) \epsilon_t \right\|_\infty \\
&= \sqrt{L} \max_{1\le j\le p}\left\| \hat \Sigma_{S_j,S_j}^{-1} \right\|_\infty \cdot \left\| \frac1n\sum_{t=1}^n\Psi_S(X_{t-1}) \epsilon_t \right\|_\infty .
\end{align*}
By Lemma \ref{lambda:l2:a}, with probability at least $1-(pL)^{-c_1}$, 
\begin{align}
\left\| \frac1n\sum_{t=1}^n\Psi_S(X_{t-1}) \epsilon_t \right\|_\infty \le c_2 \sqrt{\frac{\log(pL)}{n} } .
\end{align}
Note that
\begin{align*}
\left\|\hat\Sigma_{S,S}^{-1} \right\|_\infty = \max_{1\le j\le p} \left\|\hat\Sigma_{S_j,S_j}^{-1} \right\|_\infty  \le \max_{1\le j\le p} \left\|\hat\Sigma_{S_j,S_j}^{-1} \right\|_2 \cdot \sqrt{s_0L} =  \sqrt{s_0L} \left\|\hat\Sigma_{S,S}^{-1} \right\|_2.
\end{align*}
Then by \eqref{eq:eig2a}, with probability at least $1-(pL)^{-c_1}$,
\begin{align} \label{eq:thm4:I1}
|\I_1|_{2,\infty} 
&\le c_2 \sqrt{L} \cdot \frac{\sqrt{s_0 L}}{\phi_{\min} } \cdot \sqrt{\frac{\log(pL)}{n} } = c_2 \phi_{\min}^{-1} \frac{L\sqrt{s_0\log(pL)}} {\sqrt{n}}. 
\end{align}

For $\I_2$, by \eqref{eq:XiPhi} and Proposition \ref{prop:bddri}, we have
\begin{align} \label{eq:thm4:I2}
|\I_2|_{2,\infty} &\le \sqrt{L} \left\| \hat \Sigma_{S,S}^{-1} \right\|_\infty  \left\| \Psi_S(X_{t-1}) \right\|_\infty  \left\| r_t \right\|_\infty \le c B^2C(2\beta-1)^{-1}\phi_{\min}^{-1}s_0^{3/2} L^{3/2-\beta} .
\end{align}

For $\I_3$ part, note that for all $(j,k)\in S$,
\begin{align*}
\frac{1}{(1+o(1))\phi_{\max}} \left\| \hat g_{jk} \right\|_2^2\leq \hat g_{jk}^\top \left( \frac1n\sum_{t=1}^n \psi_{k} (X_{t-1}^{(k)}) \psi_{k} (X_{t-1}^{(k)})^\top \right)^{-1} \hat g_{jk} \leq 1.
\end{align*}
It follows that
\begin{align} \label{eq:gjk}
\left\| \hat g_{S} \right\|_\infty = \max_{(j,k)\in S} \left\| \hat g_{jk} \right\|_\infty \le \max_{(j,k)\in S} \left\| \hat g_{jk} \right\|_2 \le \sqrt{(1+o(1))\phi_{\max}}.
\end{align}
Therefore we obtain
\begin{align} \label{eq:thm4:I3}
|\I_3|_{2,\infty} &\le \frac12 \lambda\sqrt{L} \left\| \hat \Sigma_{S,S}^{-1} \right\|_\infty  \left\| \hat g_{S} \right\|_\infty \le \frac{\sqrt{(1+o(1))\phi_{\max}} } {2\phi_{\min} }\cdot \lambda \sqrt{s_0} L .
\end{align}

Combining \eqref{eq:thm4:I1}, \eqref{eq:thm4:I2} and \eqref{eq:thm4:I3}, we have, with probability at least $1-(pL)^{-c_1}$,
\begin{align}
|\hat b_S-b_S^*|_{2,\infty} &= \max _{(j,k)\in S}\|\hat b_{j,k}-b_{j,k}^*\|_{2}  \notag \\
&\le c_2 \phi_{\min}^{-1} \frac{L\sqrt{s_0\log(pL)}} {\sqrt{n}} + c B^2C (2\beta-1)^{-1}\phi_{\min}^{-1}s_0^{3/2} L^{3/2-\beta} + \frac{\sqrt{(1+o(1))\phi_{\max}} } {2\phi_{\min} }\cdot \lambda \sqrt{s_0} L .
\end{align}
By \eqref{eq:thml2:sparsity:np} and \eqref{eq:thml2:sparsity:lambda}, it follows that, on an event $\mathcal Z_1$ with probability approaching 1,
\begin{align*}
\max _{(j,k)\in S}\|\hat b_{j,k}-b_{j,k}^*\|_{2} \rightarrow 0.
\end{align*}
Since $\max _{(j,k)\in S}\| b_{j,k}^*\|_{2} > 0$ and will not converge to 0 asymptotically, \eqref{eq:signal} holds on an event $\mathcal Z_1$ with probability approaching 1.

\begin{enumerate}
\item[(ii).] Proof of \eqref{eq:subgradient}.
\end{enumerate}
Since $\hat b_{S^c}=b_{S^c}^*=0$, for all $(j,k)\in S^c$, \eqref{eq:eq} reduces to
\begin{equation*} 
\frac2n \sum_{t=1}^n \psi_{k}  (X_{t-1}^{(k)}) (\Psi_{S_j} (X_{t-1})^\top \hat b_{S_j} - X_{t}^{(j)}) +\lambda \hat g_{jk} =0.
\end{equation*}
It implies that
\begin{align*}
\hat g_{jk}=\frac{2}{\lambda} \left( \frac1n \sum_{t=1}^n \psi_{k}  (X_{t-1}^{(k)}) (\Psi_{S_j} (X_{t-1})^\top (b_{S_j}^*-\hat b_{S_j}) +  \frac1n\sum_{t=1}^n \psi_{k}  (\epsilon_{t}^{(j)}+r_{t}^{(j)})\right).
\end{align*}
By \eqref{eq:b:decom}, we have
\begin{align*}
\hat g_{jk} &= \left( \frac1n \sum_{t=1}^n \psi_{k}  (X_{t-1}^{(k)}) \Psi_{S_j} (X_{t-1})^\top \hat\Sigma_{S_j,S_j}^{-1} \right) \hat g_{S_j} \\
&\quad-\frac{2}{\lambda} \left(\frac1n \sum_{t=1}^n \psi_{k}  (X_{t-1}^{(k)}) \Psi_{S_j} (X_{t-1})^\top \hat\Sigma_{S_j,S_j}^{-1}\right) \frac1n \sum_{t=1}^n \Psi_{S_j} (X_{t-1}) \epsilon_{t}^{(j)} \\
&\quad-\frac{2}{\lambda} \left(\frac1n \sum_{t=1}^n \psi_{k}  (X_{t-1}^{(k)}) \Psi_{S_j} (X_{t-1})^\top \hat\Sigma_{S_j,S_j}^{-1}\right) \frac1n \sum_{t=1}^n \Psi_{S_j} (X_{t-1}) r_{t}^{(j)} \\
&\quad + \frac{2}{\lambda}\cdot \frac1n\sum_{t=1}^n \psi_{k}  \epsilon_{t}^{(j)} + \frac{2}{\lambda}\cdot \frac1n\sum_{t=1}^n \psi_{k}  r_{t}^{(j)} \\
&:= \II_1-\II_2-\II_3+\II_4+\II_5.
\end{align*}
Since for all $(j,k)\in S^c$,
\begin{align*}
\hat g_{jk}^\top \left( \frac1n\sum_{t=1}^n \psi_{k} (X_{t-1}^{(k)}) \psi_{k} (X_{t-1}^{(k)})^\top \right)^{-1} \hat g_{jk} \le \frac{1}{\phi_{\min}} \left\| \hat g_{jk} \right\|_2^2.
\end{align*}
It suffices to show $\max_{(j,k)\in S^c} \|\hat g_{jk}\|_2 < \sqrt{(1+o(1))\phi_{\min}}$. We now proceed to bound $\II_1,\II_2,\II_3,\II_4$ and $\II_5$.

For $\II_1$, by \eqref{eq:eig2c} and \eqref{eq:gjk},
\begin{align} 
\|\II_1\|_2 &\le \left\| \frac1n \sum_{t=1}^n \psi_{k}  (X_{t-1}^{(k)}) \Psi_{S_j} (X_{t-1})^\top \hat\Sigma_{S_j,S_j}^{-1} \right\|_2 \|\hat g_{S_j}\|_2 \notag\\
&\le (1+o(1)) \sqrt{\frac{\phi_{\min}}{\phi_{\max} } }\cdot \frac{1-\delta}{\sqrt{s_0}} \cdot \sqrt{s_0}\sqrt{\phi_{\max}} \notag\\
&\le (1+o(1)) (1-\delta)\sqrt{\phi_{\min}}. \label{eq:thm4:II1}
\end{align}

For $\II_2$, by Lemma \ref{lambda:l2:a}, as $s_0<n$, with probability at least $1-(nL)^{-c_3}$
\begin{align} \label{eq:thm4:II2}
\|\II_2\|_{2} 
&\le \frac{2}{\lambda} \cdot \sqrt{\frac{\phi_{\min}}{\phi_{\max} } }\cdot \frac{1-\delta}{\sqrt{s_0}} \cdot \sqrt{s_0L} \left\| \frac1n\sum_{t=1}^n\Psi_{S_j}(X_{t-1}) \epsilon_{t}^{(j)} \right\|_\infty \notag \\
&\le \frac{2}{\lambda} \cdot \sqrt{\frac{\phi_{\min}}{\phi_{\max} } }\cdot \frac{1-\delta}{\sqrt{s_0}} \cdot \sqrt{s_0L} \cdot
c_4 \sqrt{\frac{\log(nL)}{n} } \notag \\
&=c_5 \frac{1}{\lambda} \sqrt{\frac{L\log(nL)}{n} }.
\end{align}

For $\II_3$, by \eqref{eq:XiPhi} and Proposition \ref{prop:bddri}, we have
\begin{align} \label{eq:thm4:II3}
\|\II_3\|_{2} &\le \frac{2}{\lambda} \cdot \sqrt{\frac{\phi_{\min}}{\phi_{\max} } }\cdot \frac{1-\delta}{\sqrt{s_0}} \cdot \sqrt{s_0L} \cdot B^2 C (2\beta-1)^{-1}s_0L^{1/2-\beta} = c_6 \frac{s_0L^{1-\beta}}{\lambda}.
\end{align}

Similarly, for $\II_4$, with probability at least $1-(nL)^{-c_7}$,
\begin{align} \label{eq:thm4:II4}
\|\II_4\|_{2} 
&\le c_8 \frac{1}{\lambda} \sqrt{\frac{L\log(nL)}{n} }.
\end{align}
For $\II_5$,
\begin{align} \label{eq:thm4:II5}
\|\II_5\|_{2} &\le 2B^2 C (2\beta-1)^{-1}\frac{s_0L^{1-\beta}}{\lambda} = c_9 \frac{s_0 L^{1-\beta}}{\lambda}.
\end{align}

In view of \eqref{eq:thm4:II1}, \eqref{eq:thm4:II2}, \eqref{eq:thm4:II3}, \eqref{eq:thm4:II4} and \eqref{eq:thm4:II5}, for all $(j,k)\in S^c$, we can obtain, with probability at least $1-(nL)^{-c_3}-(nL)^{-c_7}$,
\begin{align}
\|\hat g_{jk}\|_2\le (1+o(1))(1-\delta)\sqrt{\phi_{\min}}+   (c_5+c_8) \frac{1}{\lambda} \sqrt{\frac{L\log(nL)}{n} } + (c_6+c_9) \frac{s_0L^{1-\beta}}{\lambda}. 
\end{align}
By \eqref{eq:thml2:sparsity:lambda}, it follows that, on an event $\mathcal Z_2$ with probability approaching 1,
\begin{align*}
\|\hat g_{jk}\|_2\le (1-\delta)\sqrt{\phi_{\min}}+o(1). 
\end{align*}
Hence, \eqref{eq:subgradient} holds on an event $\mathcal Z_2$ with probability approaching 1. Then Theorem \ref{thm:l2:sparsity} follows.

\end{proof}

\bibliographystyle{apalike}
\bibliography{ref}

\end{document}